\documentclass[a4paper,11pt]{article}

\usepackage[utf8]{inputenc}
\usepackage[T1]{fontenc}
\usepackage[leqno]{amsmath}
\usepackage{mathtools}
\usepackage{ifthen}
\mathtoolsset{centercolon,mathic}
\usepackage{amsthm,thmtools,amssymb,stmaryrd}
\usepackage{enumerate}
\usepackage{url}
\usepackage{mathrsfs}
\usepackage[english]{babel}
\usepackage{tikz}
\usepackage{bussproofs}
\usetikzlibrary{shapes.misc}
\usepackage{longtable}

\newtheorem{theorem}{Theorem}
\newtheorem{lemma}[theorem]{Lemma}
\newtheorem{definition}[theorem]{Definition}
\newtheorem{corollary}[theorem]{Corollary}
\newtheorem{remark}[theorem]{Remark}

\newcommand{\llor}{\,\lor\,}

\newcommand{\lland}{\,\land\,}

\newcommand{\tto}{\,\to\,}
\newcommand{\lra}{\leftrightarrow}
\newcommand{\llra}{\,\lra\,}
\newcommand{\impl}{\quad\Longrightarrow\quad}
\newcommand{\gdw}{\quad\Longleftrightarrow\quad}
\newcommand{\und}{\;\;\&\;\;}
\newcommand{\vvdash}{\,\vdash\,}

\newcommand{\Nbb}{\mathbb{N}}
\newcommand{\Sbb}{\mathbb{S}}

\newcommand{\Bcal}{\mathcal{B}}
\newcommand{\Ccal}{\mathcal{C}}
\newcommand{\Dcal}{\mathcal{D}}

\newcommand{\Fcal}{\mathcal{F}}
\newcommand{\Hcal}{\mathcal{H}}
\newcommand{\Kcal}{\mathcal{K}}
\newcommand{\Lcal}{\mathcal{L}}
\newcommand{\Pcal}{\mathcal{P}}
\newcommand{\Scal}{\mathcal{S}}

\newcommand{\ffrak}{\mathfrak{f}}
\newcommand{\mfrak}{\mathfrak{m}}
\newcommand{\nfrak}{\mathfrak{n}}

\newcommand{\BC}{\mathsf{(BC)}}

\newcommand{\varnothingunder}{\underline{\varnothing}}
\newcommand{\omegaunder}{\underline{\omega}}

\newcommand{\alphah}{{\widehat{\alpha}}}

\newcommand{\alphaxih}{{\widehat{\alpha_\xi}}}

\newcommand{\alphazh}{{\widehat{\alpha_0}}}
\newcommand{\betah}{{\widehat{\beta}}}

\newcommand{\ddfrac}[2]{\dfrac{\phantom{a}#1\phantom{a}}
                              {\phantom{a}\displaystyle{#2}\phantom{a}}}

\newcommand{\ls}{\Lcal^\star}

\newcommand{\ltset}{\Lcal^{set}_2}
\newcommand{\ltwo}{\mathrm{L}_2}

\newcommand{\pooca}{\Pi^1_1\mbox{-}\mathsf{CA}}
\newcommand{\poocas}{(\pooca^*)}

\newcommand{\cut}{(\mathsf{Cut})}
\newcommand{\sref}{(\Scal_0\mbox{-}\mathsf{Ref})}

\newcommand{\kp}{\mathsf{KP}}
\newcommand{\ido}{\mathsf{ID}_1}

\newcommand{\frs}{\mathsf{RS}}
\newcommand{\rs}{\frs^*}

\newcommand{\atros}{(\atro^*)}
\newcommand{\bis}{(\bi^*)}

\newcommand{\fpo}{\mathsf{FP}_0}

\newcommand{\BH}{\mathcal{BH}}

\newcommand{\On}{\mathit{On}}
\newcommand{\Pow}{\mathit{Pow}}

\newcommand{\vecc}[1]{\vec{#1}\hspace{0.35ex}}
\newcommand{\rk}{\mathit{rk}}
\newcommand{\lev}{\mathit{lev}}

\newcommand{\Ord}{\mathit{Ord}}
\newcommand{\Tran}{\mathit{Tran}}
\newcommand{\Succ}{\mathit{Succ}}
\newcommand{\FinOrd}{\mathit{FinOrd}}

\newcommand{\Ma}{M_\alpha}
\newcommand{\Mb}{M_\beta}

\newcommand{\natsum}{\mathrel{\#}}

\def\prov#1#2{\mathrel{\vrule width 0.3pt height 8pt depth
2pt\mkern-2mu
\textstyle{\frac{\hskip0.5ex\raise2pt\hbox{$\scriptstyle#1$}\hskip0.5ex}
{\hskip0.5ex#2\hskip0.5ex}}}}

\newcommand{\provv}[2]{\prov{#1}{    \stackrel{}{#2}           }}

\newcommand{\proofrule}[2]
{
\displaystyle{\frac{#1}{#2}}
}

\newcommand{\TT}{\pooca + \mathit{TI}({<}\BH)}
\newcommand{\TTa}{\pooca + \mathit{TI}(\alpha{+}\omega)}

\newcommand{\atro}{\mathsf{ATR}_0}
\newcommand{\atroset}{\mathsf{ATR}^{set}_0}

\newcommand{\bi}{\mathsf{BI}}

\newcommand{\truth}{\mathcal{T}}

\newbox\gnBoxA
\newdimen\gnCornerHgt
\setbox\gnBoxA=\hbox{$\ulcorner$}
\global\gnCornerHgt=\ht\gnBoxA
\newdimen\gnArgHgt
\def\Godelnum #1{%
\setbox\gnBoxA=\hbox{$#1$}%
\gnArgHgt=\ht\gnBoxA%
\ifnum     \gnArgHgt<\gnCornerHgt \gnArgHgt=0pt%
\else \advance \gnArgHgt by -\gnCornerHgt%
\fi \raise\gnArgHgt\hbox{$\ulcorner$} \box\gnBoxA %
\raise\gnArgHgt\hbox{$\urcorner$}}

\newcommand{\goed}[1]{\Godelnum{ #1}}

\title{Admissible extensions of subtheories of\\ second order arithmetic}

\author{Gerhard J\"ager and Michael Rathjen}

\begin{document}

\maketitle

\begin{abstract}\noindent
In this paper we study admissible extensions of several theories $T$ of reverse mathematics. The idea is that in such an extension the structure $\mathfrak{M} =(\Nbb,\Sbb,\in)$ of the natural numbers and collection of sets of natural numbers $\Sbb$ has to obey the axioms of $T$ while simultaneously one also has a set-theoretic world with transfinite levels erected on top of
$\mathfrak{M}$ governed by the axioms of Kripke-Platek set theory, $\kp$. 

In some respects, the admissible extension of $T$ can be viewed as a proof-theoretic analog of Barwise's admissible cover of an arbitrary model of set theory; see \cite{barwise75}. However, by contrast, the admissible extension of $T$ is usually not a conservative extension of $T$. Owing to the interplay of $T$ and $\kp$, either theory's axioms may force new sets of natural to exists which in turn may engender yet new sets of naturals on account of the axioms of the other.

The paper discerns a general pattern though. It turns out that for many familiar theories $T$, the second order part of the admissible cover of $T$ equates to $T$ augmented by transfinite induction over all initial segments of the Bachmann-Howard ordinal. Technically, the paper uses a novel type of ordinal analysis, expanding that for $\kp$ to the higher set-theoretic universe while at the same time treating the world of subsets of $\Nbb$ as an unanalyzed class-sized urelement structure.

Among the systems of reverse mathematics, for which we determine the admissible extension, are $\Pi^1_1\mbox{-}\mathsf{CA}_0$ and $\mathsf{ATR}_0$ as well as the theory of bar induction, $\mathsf{BI}$.

\par\medskip


\par\medskip

\end{abstract}

\section{Introduction}

In this paper we study admissible extensions of several theories $T$ of reverse mathematics. The idea is that in such an extension the structure $\mathfrak{M} =(\Nbb,\Sbb,\in)$ of the natural numbers and collection of sets of natural numbers $\Sbb$ has to obey the axioms of $T$ while simultaneously one also has a set-theoretic world with transfinite levels erected on top of
$\mathfrak{M}$ governed by the axioms of Kripke-Platek set theory, $\kp$. 

In some respects, the admissible extension of $T$ can be viewed as a proof-theoretic analog of Barwise's admissible cover of an arbitrary model of set theory (in his book ``Admissible Sets and Structures''). However, by contrast, the admissible extension of $T$ is usually not a conservative extension of $T$. Owing to the interplay of $T$ and $\kp$, either theory's axioms may force new sets of naturals to exists which in turn may engender yet new sets of naturals on account of the axioms of the other.

This approach will be studied in detail and paradigmatically by combining $\Pi^1_1$ comprehension on the natural numbers with Kripke-Platek set theory, the respective theory being called $\kp +\poocas$. In the next two sections we present the syntactic machinery of this system and then turn to its Tait-style reformulation, which is convenient for the later proof-theoretic analysis. For this proof-theoretic analysis we make use of a more or less standard ordinal notation system, following Buchholz \cite{buchholz86a}, that includes the Bachmann-Howard ordinal $\BH$ and the collapsing functions needed later.

Section~\ref{s:rs} introduces the semi-formal system $\rs$ in which $\kp +\poocas$ will be embedded. For the subsequent analysis of $\rs$ the paper uses a novel type of ordinal analysis, expanding that for $\kp$ to the higher set-theoretic universe while at the same time treating the world of subsets of $\Nbb$ as an unanalyzed class-sized urelement structure.

Based on these results we turn to the reduction of $\kp +\poocas$ to $\TT$ in Section~\ref{s:reduction}. Main ingredients here are the so-called $\alpha$-trees - extending Simpson's suitable trees in \cite{simpson09} - and a specific truth definition based on those. 

The final section shows how the previous results can be extended to other systems of reverse mathematics such as $\atro$ and $\bi$.

\section{$\kp$ and $\poocas$}

Let $\Lcal_\in$ be the usual language of set theory with $\in$ as its only non-logical relation symbol plus set constants $\varnothingunder$ (for the empty set) and $\omegaunder$ (for the first infinite ordinal). $\ltset$ is the extension of $\Lcal_\in$ by countably many unary relation variables.
As we will see, $\ltset$ may be regarded as a two-sorted first order language where the first sort is supposed to range over sets and the second sort over collections of natural numbers. 

The \emph{set terms} of $\ltset$ are the set constants and the set variables. The \emph{atomic formulas} of $\ltset$ are the expressions of the form $(a \in b)$, $(a \notin b)$, $U(a)$, and $\neg U(a)$, where $a,b$ are set terms and $U$ is a relation variable. The \emph{formulas} of $\ltset$ are built up from these atomic formulas by means of the propositional connectives $\lor, \land$, the bounded set quantifiers $(\exists x \in a)$, $(\forall x \in a)$, the unbounded set quantifiers $\exists x$, $\forall x$, and the relation quantifiers $\exists X$, $\forall X$. We use as metavariables (possibly with subscripts):
\begin{itemize}
\item
$u,v,w,x,y,z$ for set variables,
\item 
$a,b,c$ for set terms,
\item
U,V,W,X,Y,Z for relation variables,
\item
$A,B,C,D$ for formulas.
\end{itemize}
Observe that all $\ltset$ formulas are in negation normal form. We assume classical logic throughout this article. Therefore, the negation $\neg A$ of an $\ltset$ formula $A$ is defined via de Morgan's laws and the law of double negation. Furthermore, $(A \to B)$ is defined as $(\neg A \lor B)$ and $(A \lra B)$ as $((A \to B) \land (B \to A))$. 

To simplify the notation we often omit parentheses if there is no danger of confusion. Equality is not a basic symbol but defined for sets and relations:
\begin{align*}
   (a=b) &\;:=\; (\forall x \in a)(x \in b) \lland (\forall x \in b)(x \in a),
   \\[1.5ex]
   (U=V) &\;:=\; (\forall x \in \omegaunder)(U(x) \lra V(x)).%
   \end{align*}
Since relations are supposed to range over subcollections of $\omegaunder$ -- see  axiom (Sub $\omegaunder)$ below -- this definition of the equality of relations makes sense.
   
The vector notation $\vec{a}$ will be used to denote finite sequences of $\ltset$ terms. $A^{(a)}$ results from $A$ by restricting all unbounded set quantifiers to $a$; relation quantifiers $QX$ are not affected by this restriction. The $\Delta_0$, $\Sigma$, $\Pi$, $\Sigma_n$, and $\Pi_n$ formulas of $\ltset$ are defined as usual where relation variables are permitted as parameters. However, relation quantifiers are not allowed in these formulas.
The \emph{bounded} formulas are all those $\ltset$ formulas which do not contain unbounded set quantifiers; they may contain relation quantifiers.

Moreover, we shall employ the common set-theoretic terminology and the standard notational conventions, for example:
\begin{itemize}
\item
$\Tran[a] \;:=\; (\forall x \in a)(\forall y \in x)(y \in a)$,
\item
$\Ord[a] \;:=\; \Tran[a] \lland (\forall x \in a)\Tran[x]$,
\item
$\Succ[a] \;:=\; \Ord[a] \lland (\exists x \in a)(a = x \cup \{x\})$, 
\item
$\FinOrd[a] \;:=\;  \Ord[a] \lland (a = \varnothingunder \llor \Succ[a]) \lland
(\forall x \in a)(x = \varnothingunder \llor \Succ[x]).$

\end{itemize} 
Clearly, $\Succ[a]$ expresses that $a$ is a successor ordinal and $\FinOrd[a]$ says that $a$ is a finite ordinal.

Now we formulate \emph{Kripke-Platek set theory} $\kp$ in $\ltset$, based on classical logic.%
\footnote{Observe that $(a=a)$ and $(U=U)$ are provable by logic.}
Its non-logical axioms are:
\[
\begin{array}{lll}
\mbox{(Equality)} &  a=b \land D[a] \;\to\; D[b].
\phantom{aaaaaaaaaaaaaaaaaaaaaaaaaaaaaa}
\\[1.5ex]
\mbox{(Sub $\omegaunder)$ } & U(a) \tto a \in \omegaunder.
\\[1.5ex]
\mbox{(Pair)} &  \exists z(a \in z \lland b \in z).
\\[1.5ex]
\mbox{(Union)} &  \exists z(\forall y \in a)(\forall x \in y)(x \in z).
\\[1.5ex]
\mbox{(Empty set)} & (\forall x \in\varnothingunder)(x \neq x).
\\[1.5ex]
\mbox{(Infinity)} &  a \in \omegaunder \;\lra\; \FinOrd[a].
\\[1.5ex]
\mbox{($\Delta_0$-Sep)} &  \exists z(z = \{x \in a : D[x]\}).
\\[1.5ex]
\mbox{($\Delta_0$-Col)} &  (\forall x \in a)\exists y D[x,y] \tto 
                        \exists z(\forall x \in a)(\exists y \in z)D[x,y].
\\[1.5ex]
\mbox{($\in$-Ind)} &  \forall x((\forall y \in x)A[y] \tto A[x]) \,\to\, \forall x A[x].
\end{array}
\]
The formulas $D$ in the schemas (Equality), ($\Delta_0$-Sep), and ($\Delta_0$-Col) are $\Delta_0$ whereas the formula $A$ in ($\in$-Ind) ranges over arbitrary formulas of $\ltset$.

\medskip

It is easy to see that the theory $\kp$ is a conservative extension of the usual first order formalization of Kripke-Platek set theory with infinity. This system is called $\kp \omega$ in some texts.

In this article we are interested in the theory that is obtained from $\kp$ by adding the following form of $\Pi^1_1$ comprehension, where $D[x,Y]$ is a $\Delta_0$ formula of $\ltset$,
\begin{gather*}
\exists Z (\forall x \in \omegaunder)(Z(x) \;\lra\; \forall Y D[x,Y]). \tag*{$\poocas$}
\end{gather*}

\medskip

From $(\Delta_0$-Sep) -- and here it is crucial that relations are permitted as parameters -- we immediately obtain that  the intersection of any relation with $\omegaunder$ is extensionally equal to a subset of $\omegaunder$. On the other hand, in view of $\poocas$ we also know that for every subset of $\omegaunder$ there exists a relation with the same elements. 

Let $\ltwo$ be the usual language of second order arithmetic as presented, for example, in Simpson \cite{simpson09}. A natural translation of $\ltwo$ into $\ltset$ is defined as follows: The number variables of $\ltwo$ are interpreted in $\ltset$ as ranging over $\omegaunder$ and the set variables of $\ltwo$ are replaced in $\ltset$ by relation variables. The other symbols and connectives of $\ltwo$ are dealt with as in \cite[VII.3.9]{simpson09}.

It is also clear that our  comprehension schema $\poocas$ implies the usual form of $\Pi^1_1$ comprehension of second order arithmetic modulo its natural translation into $\ltset$.

\begin{remark}\rm
Let $A'[x,Y]$ be the natural translation of an arithmetical $A[x,Y]$. Then $\kp +\poocas$ proves that
\[ (\exists z \subseteq \omegaunder)(\forall x \in \omegaunder)(x \in z \;\lra\;
   (\forall y \subset \omegaunder)A'[x,y]). \]
\end{remark}

\section{A Tait-style reformulation of $\kp +\poocas$}

The basic idea of a Tait-style calculus is that it derives finite sets of $\ltset$ formulas rather than individual formulas. As in one-sided sequent calculi, the intended meaning of such a set of formulas is the disjunction of its elements. The following reformulation of $\kp +\poocas$ is for technical reasons only and simplifies the proof of Theorem~\ref{t:embedding}.

The capital Greek letters $\Gamma,\Theta,\Lambda$ (possibly with subscripts) will act as metavariables for finite sets of $\ltset$ formulas. Also, we write (for example) $\Gamma, A_1, \ldots, A_n$ for the set $\Gamma \cup \{A_1,\ldots,A_n\}$; similarly for expressions such as $\Gamma,\Theta,A$.

\bigskip\medskip

\noindent\textbf{The Tait-style axioms of $\kp + \poocas$}
\[ \begin{array}{ll}
\mbox{(TnD)} & \Gamma,\, \neg D,\, D.
\phantom{aaaaaaaaaaaaaaaaaa}
\\[1.5ex]
\mbox{(Equality)} & \Gamma,\, a \neq b,\, \neg D[a],\, D[b]
\\[1.5ex]
\mbox{(Sub $\omegaunder)$ } & \Gamma,\, \neg U(a),\, a \in \omegaunder.
\\[1.5ex]
\mbox{(Pair)} & \Gamma,\, \exists z(a \in z \lland b\in z).
\\[1.5ex]
\mbox{(Union)} & \Gamma,\, \exists z(\forall y \in a)(\forall x \in y)(x \in z).
\\[1.5ex]
\mbox{(Empty set)} & \Gamma,\, a  \not\in \varnothingunder.
\\[1.5ex]
\mbox{(Infinity)} & \Gamma,\, a \in \omegaunder \llra \FinOrd[a].
\\[1.5ex]
\mbox{($\Delta_0$-Sep)} & \Gamma,\, \exists z(z = \{x \in a : D[x]\}).
\\[1.5ex]
\mbox{($\Delta_0$-Col)} & \Gamma,\, (\forall x \in a)\exists y D[x,y] \tto 
                        \exists z(\forall x \in a)(\exists y \in z)D[x,y].
\\[1.5ex]
\mbox{($\in$-Ind)} & \Gamma,\, \forall x((\forall y \in x)A[y] \tto A[x]) \,\to\, \forall x A[x].
\\[1.5ex]
\poocas & \Gamma,\, \exists Z (\forall x \in \omegaunder)(Z(x) \;\lra\; \forall Y D[x,Y]).
\end{array}
\]
In these axioms the formulas $D$ in (TnD), (Equality), ($\Delta_0$-Sep), ($\Delta_0$-Col), and $\poocas$ are supposed to be $\Delta_0$. The formula $A$ in ($\in$-Ind) may be an arbitrary $\ltset$ formula. (TnD) stands for ``Tertium non datur''.

\medskip\bigskip

\noindent\textbf{ The Tait-style Inference rules of $\kp +\poocas$}
\medskip

\allowdisplaybreaks

\[
\begin{array}{ccccc}
\ddfrac{\Gamma,\, A,\, B}{\Gamma,\, A \lor B} & (\lor)  & &
\ddfrac{\Gamma,\, A \qquad \Gamma,\, \, B}{\Gamma,\, A \land B} &(\land)
\\[3.5ex]
\ddfrac{\Gamma,\, A[b]}{\Gamma,\, \exists x A[x]}\ &(\exists) & &
\ddfrac{\Gamma,\, A[u]}{\Gamma,\, \forall x A[x]} & (\forall)
\\[3.5ex]
\ddfrac{\Gamma,\, b \in a \lland A[b]}{\Gamma,\, (\exists x \in a) A[x]} & (b\exists) & &
\ddfrac{\Gamma,\, u \in a  \tto A[u]}{\Gamma,\, (\forall x \in a)A[x]} & (b\forall)
\\[3.5ex]
\ddfrac{\Gamma,\, A[V]}{\Gamma,\, \exists X A[X]} & (\exists_2) & &
\ddfrac{\Gamma,\, A[U]}{\Gamma,\, \forall X A[X]} & (\forall_2)
\\[3.5ex]
\ddfrac{\Gamma,\, A \qquad \Gamma,\, \neg A}{\Gamma} & \cut & &
\phantom{a} & \phantom{a}
\end{array} \]

\smallskip

\noindent Of course, it is demanded that in $(\forall)$ and $(b\forall)$ the eigenvariable $u$ must not occur in the conclusion; the same is the case for the variable $U$ in $(\forall_2)$.

We say that $\Gamma$ is Tait-style devivable from $\kp + \poocas$ iff there exists a finite sequence of finite sets of $\ltset$ formulas
\[ \Theta_0,\ldots,\Theta_h  \]
such that $\Theta_h$ is the set $\Gamma$ and for any $i = 0,\ldots, h$ one of the following two conditions is satisfied:
\begin{itemize}
\item 
$\Theta_i$ is a Tait-style axiom of $\kp + \poocas$.
\item $\Theta_i$ is the conclusion of an inference of a Tait-style inference rule of  $\kp + \poocas$ with premise(s) from $\Theta_0,\ldots,\Theta_{i-1}$.
\end{itemize}
In this case we write $\kp + \poocas \vdash^h \Gamma$ and say that $\Gamma$ has a proof of length $h$. It is an easy exercise to show that a formula $A$ is provable in one of the usual Hilbert-style formalizations of $\kp +\poocas$ iff $\kp + \poocas \vdash^h A$ for some natural number $h$. Details are left to the reader.

\section{Ordinal notations}

In the next sections we establish the upper proof-theoretic bound of the theory $\kp +\poocas$.
Our method of choice is the ordinal analysis of $\kp +\poocas$ via a system $\rs$ of ramified set theory. And in order to build up this system and to control the derivations in $\rs$ we work with specific ordinal notations. 

Buchholz has developed several ordinal notation systems based on so-called collapsing functions; see, for example, Buchholz \cite{buchholz86a,buchholz92a,buchholz02a}. In the following we will work with a reduced version, which is sufficient for our purposes. 

Let $\On$ be the collection of all ordinals and let $\Omega$ be a sufficiently large ordinal. To simplify matters we set $\Omega := \aleph_1$, but also $\omega_1^{ck}$ or even somewhat small ordinals could do the job. The following outline is based on \cite{buchholz92a}.

\begin{definition}\rm
The set of ordinals $C(\alpha,\beta)$ and the ordinals $\psi(\alpha)$ are defined for all ordinals $\alpha$ and $\beta$ by induction on $\alpha$.
\begin{enumerate}[(i)]
\item
$\{0,\Omega\} \cup \beta \subseteq C(\alpha,\beta)$.
\item
If $\eta,\xi \in C(\alpha,\beta)$, then $\eta +\xi \in C(\alpha,\beta)$ and $\omega^\xi \in C(\alpha,\beta)$.
\item
If $\xi <\alpha$ and $\xi \in C(\alpha,\beta)$, then $\psi(\xi) \in C(\alpha,\beta)$.
\item
$\psi(\alpha) \,:=\, \min(\{ \eta \in \On : C(\alpha,\eta) \cap \Omega = \eta \})$.
\end{enumerate}
\end{definition}

The following lemma summarizes some key properties of the sets $C(\alpha,\beta)$ and the function $\psi$. For its proof see Buchholz \cite{buchholz86a}.

\begin{lemma}  \label{l:o0}
We have for all ordinals $\alpha,\alpha_1,\alpha_2,\beta$:
\begin{enumerate}[(1)]
\item
$\psi(\alpha) < \Omega$.
\item
$C(\alpha,\psi(\alpha)) \cap \Omega = \psi(\alpha)$.
\item
$\psi(\alpha)$ is an $\varepsilon$-number.
\item
If $\alpha_1 <\alpha_2$ and $\alpha_1 \in C(\alpha_2,\psi(\alpha_2))$, then $\psi(\alpha_1) < \psi(\alpha_2)$.
\item
If $\alpha_1 \leq \alpha_2$, then $\psi(\alpha_1) \leq \psi(\alpha_2)$ and
$C(\alpha_1,\psi(\alpha_1)) \subseteq C(\alpha_2,\psi(\alpha_2))$.
\item
$C(\alpha,0) = C(\alpha,\psi(\alpha))$.

\end{enumerate}
\end{lemma}

We write $\varepsilon_{\Omega+1}$ for the least ordinal $\alpha > \Omega$ such that $\omega^\alpha = \alpha$. Its collapse $\BH := \psi(\varepsilon_{\Omega+1})$ is called the \emph{Bachmann-Howard ordinal}. This number gained importance in proof theory since it is the proof-theoretic ordinal of the theory $\ido$ of one positive inductive definition and of Kripke-Platek set theory $\kp$; see, for example, Buchholz and Pohlers \cite{buchholz-pohlers78a}, J\"ager \cite{j82a}, and Pohlers \cite{pohlers81a}.

\section{The semi-formal ramified system $\rs$} \label{s:rs}

In this section we introduce the semi-formal proof system $\rs$ of ramified set theory. We begin with extending our language $\ltset$ to the language $\ls$ and then present the axioms and rules of inference of $\rs$. Afterwards we turn to operator controlled derivations and some basic properties of $\rs$. We show that $\kp + \poocas$ can be embedded into $\rs$ and prove cut elimination and collapsing for $\rs$. Henceforth, all ordinals used in this section on the metalevel range over the set $C(\varepsilon_{\Omega+1},0)$ if not stated otherwise.

\subsection{The language $\ls$}

The basic idea is to extend the language $\ltset$ to the language $\ls$ by adding unary relation symbols $\Ma$ and new quantifiers $\exists x^\alpha\, /\,  \forall x^\alpha$ for all $\alpha < \Omega$. The quantifiers $Q x^\alpha$ are supposed to range over $\Ma$, and later, see Subsection~\ref{s:suitable}, an atomic formula $\Ma(a)$ will be interpreted as stating that $a$ can be coded by a so-called $\alpha$-tree, see Definition~\ref{MR7} below.
\begin{definition}\rm \label{d:rs-formulas}
The \emph{formulas} $F$ of $\rs$, their \emph{ranks} $\rk(F)$ and \emph{parameter sets} $|F|$ are inductively defined as follows.
\begin{enumerate}
\item
If $a$ and $b$ are set terms of $\ltset$, then $(a\in b)$ and $(a \notin b)$ are formulas with
\[ \rk(a \in b) := \rk(a \notin b) := 0\; \mbox{and}\;  |(a \in b)| := |(a \notin b)| := \{0\}. \]
\item
If $a$ is a set term and $U$ a relation variable, then $U(a)$ and $\neg U(a)$ are formulas with
\[ rk(U(a)) := \rk(\neg U(a)) := 0\; \mbox{and}\; |U(a)| := |\neg U(a)| := \{0\}. \]
\item
If $a$ is a set term and $\Ma$ one of these new relation symbols, then $\Ma(a)$ and $\neg \Ma(a)$ are formulas with
\[ rk(\Ma(a)) := \rk(\neg \Ma(a)) := \omega\alpha\; \mbox{and}\; 
   |\Ma(a)| := |\neg \Ma(a)| := \{\alpha\}. \]
\item
If $F$ and $G$ are formulas, then $(F \lor G)$ and $(F \land G)$ are formulas with 
\begin{gather*}
rk(F \lor G) := \rk(F \land G) := \max(\rk(F),\rk(G))+1,
\\[1ex]
|F \lor G| := |F \land G| := |F| \cup |G|. 
\end{gather*}
\item
If $F[u]$ is a formula, then $(\exists x \in a)F[x]$ and $(\forall x \in a)F[x]$ are formulas with
\begin{gather*}
\rk((\exists x  \in a)F[x]) :=  \rk((\forall x  \in a)F[x]) := \rk(u \in a \lland F[u]) +1.
\\[1ex]
|(\exists x  \in a)F[x]| := |(\forall x  \in a)F[x]| := |F[u]|.
\end{gather*}
\item
If $F[u]$ is a formula, then $\exists x^\alpha F[x]$ and $\forall x^\alpha F[x]$ are formulas with 
\begin{gather*}
\rk(\exists x^\alpha F[x]) :=  \rk(\forall x^\alpha F[x]) := \rk(\Ma(u) \lland F[u]) +1.
\\[1ex]
|\exists x^\alpha F[x]| := |\forall x^\alpha F[x]| := \{\alpha\} \cup |F[u]|.
\end{gather*}
\item
If $F[u]$ is a $\Delta_0$ formula of $\ltset$, then $\exists x F[x]$ and $\forall x F[x]$ are formulas with
\begin{gather*}
\rk(\exists x F[x]) :=  \rk(\forall x F[x]) := \Omega, 
\\[1ex]
|\exists x F[x]| := |\forall x F[x]| := |F[u]|. 
\end{gather*}
\item
If $F[u]$ is not a $\Delta_0$ formula of $\ltset$, then $\exists x F[x]$ and $\forall x F[x]$ are also formulas but with
\begin{gather*}
\rk(\exists x F[x]) :=  \rk(\forall x F[x]) := \max(\Omega+1,\rk(F[u])+3), 
\\[1ex]
|\exists x F[x]| := |\forall x F[x]| := |F[u]|. 
\end{gather*}
\item
If $ F[U]$ is a formula, then $\exists X F[X]$ and $\forall X F[X]$ are formulas with 
\begin{gather*}
\rk(\exists X F[X]) :=  \rk(\forall X F[X]) := \rk(F[U])+1,
\\[1ex]
|\exists X F[X]| := |\forall X F[X]| := |F[U]|.
\end{gather*}
\end{enumerate}
Finally, we define the \emph{level} of a formula as 
\[
\lev(F) \;:=\; \left\{
\begin{array}{ll}
\max(|F|) & \mbox{if $\rk(F) < \Omega$},
\\[1ex]
\Omega & \mbox{if $\Omega \leq \rk(F)$}.
\end{array}
\right. \]
\end{definition}
To be precise: If $F$ is a formula of $\rs$, then $|F|$ collects the levels of the relations symbols $\Ma$ and of the quantifiers $Q x^\alpha$ occurring in $F$ plus possibly the number $0$.

There are several collections of $\ls$ formulas that will play an important role later.
\begin{enumerate}
\item
$\Dcal$ is the collection of all $\ls$ formulas in which unbounded set quantifiers $Q x$ and relation quantifiers $Q X$ do not occur. 
\item
$\Scal$ is the closure of $\Dcal$ under the propositional connectives, quantifiers $Q x^\alpha$, bounded quantifiers $(Qx \in r)$,  and unbounded existential set quantifiers. 
\item
$\Scal_0$ is the subclass of $\Scal$ that contains all $\Sigma$ formulas of $\ltset$ that are not $\Delta_0$.
\item
$\Bcal$ consists of all $\ls$ formulas that do not contain unbounded set quantifiers $Q x$. \end{enumerate}
Some important properties of the ranks of $\rs$ formulas are summarized in the following lemma. Its proof is straightforward and will be omitted.

\begin{lemma}\rm\label{l:rank} \quad
\begin{enumerate}[(1)]
\item
$\rk(F) \,=\, \rk(\neg F) \,<\, \omega \cdot \lev(F) + \omega \,\leq\, \Omega+\omega$.
\item
$\rk(F) < \Omega$ iff $F$ belongs to $\Bcal$.
\item
If $\rk(\exists x F[x]) = \Omega$, then $F[u]$ is a $\Delta_0$ fromula of $\ltset$.
\item
$\rk(\Ma(a) \lland F[a]) \,<\, \rk(\exists x^\alpha F[x]) \,<\, \rk(\exists x F[x])$.
\item
$\rk(s \in r \lland F[s]) \,<\, \rk((\exists x \in r)F[x]) \,<\, \rk(\exists x F[x])$.
\item
$\rk(F[P]) \,=\, \rk(F[U]) \,<\, \rk(\exists X F[X])$.

\end{enumerate}
\end{lemma}

\subsection{Axioms and rules of inference of $\rs$}

We shall use a Tait-style calculus as proof system for $\rs$. In the following the Greek capital letters $\Gamma,\Theta,\Lambda$ (possibly with subscripts) will now act as metavariables for finite sets of $\ls$ formulas. If  $\Gamma$ is the set $\{F_1,\ldots,F_k\}$ we define
\[ |\Gamma| \;:=\; |F_1| \cup \ldots \cup |F_k| \quad\mbox{and}\quad
   \Gamma^\lor \;:=\; F_1 \lor \ldots \lor F_k  \]
such that $|\Gamma|$ collects the parameter sets of the formulas in $\Gamma$ and $\Gamma^\lor$ is the disjunction of its elements. 
\bigskip

\noindent\textbf{Axioms of $\rs$}

\begin{enumerate}[(1)]

\vspace{0.5ex}\item\label{a:1}
$\Gamma,\, \neg F,\, F$\quad for $F \in \Bcal$.

\vspace{0.5ex}\item
$\Gamma,\, a \neq b,\, \neg F[a],\, F[b]$\quad for $F[u] \in \Bcal$.

\vspace{0.5ex}\item
$\Gamma,\, \neg U(a),\, a \in \omegaunder$.

\vspace{0.5ex}\item\label{a:5}
$\Gamma,\, \neg M_0(a)$.

\vspace{0.5ex}\item \label{a:6}
$\Gamma,\, M_1(\varnothingunder)$.

\vspace{0.5ex}\item
$\Gamma,\, a \notin \varnothing$.

\vspace{1ex}\item \label{a:8}
$\Gamma,\, M_{\omega{+}1}(\omegaunder)$.

\vspace{1ex}\item
$\Gamma,\, (a \in \omegaunder \;\lra\; \FinOrd[a])$.

\vspace{1ex}\item
$\Gamma,\, \neg \Ma(a),\, \Mb(a)$\quad for $\alpha \leq \beta$.

\vspace{1ex}\item\label{a:12}
$\Gamma,\, \neg M_{\alpha{+}1}(a),\, b \notin a,\, \Ma(b)$.

\vspace{1ex}\item\label{a:13}
$\Gamma,\, \neg \Ma(a),\, \neg\Mb(b),\, \exists z^{\beta{+}1}(a \in z \lland b  \in z)$\quad
for $\alpha \leq \beta$.

\vspace{1ex}\item
$\Gamma,\, \neg \Ma(a),\, \exists z^\alpha(\forall y \in a)(\forall x \in y)(x \in z)$.

\vspace{1ex}\item\label{a:15}
$\Gamma,\, \neg \Ma(a),\, \exists z^{\alpha{+}1}(z = \{x \in a: D[x]\})$
\quad for $D[u]$ from  $\Delta_0$ of $\ltset$.

\vspace{1ex}\item\label{a:16}
$\Gamma,\, \exists Z(\forall x \in \omegaunder)(Z(x) \llra  \forall Y D[x,Y])$
\quad for $D[u,V]$ from  $\Delta_0$ of $\ltset$.

\end{enumerate}

Please observe that the main formulas of all axioms belong to $\Bcal$. This will be important in the later subsections when it comes to the ordinal analysis of the derivations in $\rs$. 

\bigskip

\noindent\textbf{Inference rules of $\rs$}

\begin{enumerate}[({RRRRR}1)]

\vspace{1ex}\item[$(\lor)$]\phantom{a}
$\proofrule{\Gamma,\,F,\, G}
{\Gamma,\, F \lor G}$

\vspace{1ex}\item[$(\land)$]\phantom{a}
$\proofrule{\Gamma,\, F \qquad \Gamma,\, G}
{\Gamma,\, F \land G}$

\vspace{1ex}\item[$(\neg M)$]\phantom{a}
$\proofrule{\Gamma,\, \neg \Mb(a) \,\,\, \mbox{for all $\beta < \lambda$}}
{\Gamma,\, \neg M_\lambda(a)}$\quad for $\lambda$ limit

\vspace{1ex}\item[$(\exists)$]\phantom{a}
$\proofrule{\Gamma,\, \Mb(a) \lland F[a]}
{\Gamma,\, \exists x F[x]}$

\vspace{1ex}\item[$(\forall)$]\phantom{a}
$\proofrule{\Gamma,\, \neg \Mb(a) \llor F[a] \,\,\, \mbox{for all $\beta < \Omega$ and all $a$}}
{\Gamma,\, \forall x F[x]}$

\vspace{1ex}\item[$(\exists^\alpha)$]\phantom{a}
$\proofrule{\Gamma,\, \Mb(a) \lland F[a]}{\Gamma,\, \exists x^\alpha F[x]}$\quad 
\mbox{for $\beta \leq \alpha$}

\vspace{1ex}\item[$(\forall^\alpha)$]\phantom{a}
$\proofrule{\Gamma,\, \neg\Mb(a) \llor F[a]\,\,\, \mbox{for all $\beta \leq \alpha$ and all $a$}}
{\Gamma,\, \forall x^\alpha F[x]}$

\vspace{1ex}\item[$(b\exists)$]\phantom{a}
$\proofrule{\Gamma,\, b \in a \lland F[b]}
{\Gamma,\, (\exists x \in a)F[x]}$

\vspace{1ex}\item[$(b\forall)$]\phantom{a}
$\proofrule{\Gamma,\, b \notin a \llor F[b]\,\,\, \mbox{for all $b$}}
{\Gamma,\, (\forall x \in a)F[x]}$

\vspace{1ex}\item[$(\exists_2)$]\phantom{a}
$\proofrule{\Gamma,\, F[U]}
{\Gamma,\, \exists X F[X]}$

\vspace{1ex}\item[$(\forall_2)$]\phantom{a}
$\proofrule{\Gamma,\, F[U] \,\,\, \mbox{for all $U$}}
{\Gamma,\, \forall X F[X]}$

\vspace{1ex}\item[$\cut$]\phantom{a}
$\proofrule{\Gamma,\, F \qquad \Gamma,\, \neg F}
{\Gamma}$

\vspace{1ex}\item[$\sref$]\phantom{a}
$\proofrule{\Gamma,\, F}
{\Gamma,\, \exists z F^{(z)}}$\quad for $F \in \Scal_0$

\vspace{1ex}\item[$\BC$]\phantom{a}
$\proofrule{\Gamma,\, F^\beta}{\Gamma,\, \exists z^{\beta{+}\omega} F^{(z)}}$\quad 
 for $F \in \Scal_0$

\end{enumerate}
If $F$ is from $\Scal$, we write $F^\beta$ for the result of replacing each unbounded existential quantifier $\exists x$ by $\exists x^\beta$.  This must not be confused with $F^{(a)}$ where each unbounded set quantifier $Qx$ in $F$ is replaced by $(Q x \in a)$.

\smallskip

The meaning of the rules $(\lor)$ -- $\sref$ should be self-explaining. Rule $\BC$ will be needed in connection with the boundedness and collapsing results in subsection 
\ref{ss:collapsing}.

\subsection{Derivation operators} \label{s:derivation-operators}

The general theory of derivation operators and operator controlled derivations has been introduced in Buchholz \cite{buchholz92a}. In the following we adapt his general approach to the more specific (and simpler) situation with which we have to deal here.

\begin{definition}\rm
Let $\Pow(\On)$ denote the collection of all sets of ordinals. A class function
\[ \Hcal : \Pow(\On) \to \Pow(\On) \]
is called a \emph{derivation operator d-operator for short)}  iff it is a closure operator and satisfies the following conditions for all $X,Y \in \Pow(\On)$:
\begin{enumerate}[(i)]
\item
$X \subseteq \Hcal(X)$.
\item
$Y \subseteq \Hcal(X) \impl \Hcal(Y) \subseteq \Hcal(X)$.
\item
$\{0,\Omega\} \subseteq \Hcal(X)$.
\item
If $\alpha$ has Cantor normal form $\omega^{\alpha_1} + \ldots + \omega^{\alpha_k}$, then
\[ \alpha \in \Hcal(X) \gdw \alpha_1,\ldots,\alpha_k \in \Hcal(X). \]
\end{enumerate}
\end{definition}

These requirements ensure that every d-operator $\Hcal$ is monotone, inclusive, and idempotent.  Every $\Hcal(X)$ is closed under $+$ as well as $\xi \mapsto \omega^\xi$, and under the decomposition of its members into their Cantor normal form components.

Let $\Hcal$ be a d-operator. Then we define for all finite sets of ordinals $\mfrak$ the operators
\[ \Hcal[\mfrak] \;:\; \Pow(\On) \to \Pow(\On) \]
by setting for all $X \in \Pow(\On)$:
\[ \Hcal[\mfrak](X) \;:=\; \Hcal(\mfrak \cup X). \]
Also, we simply write $\Hcal[\sigma]$ for $\Hcal[\{\sigma\}]$. If $\Hcal$ and $\Kcal$ are d-operators, then we write $\Hcal \subseteq \Kcal$ to state that
\[ \Hcal(X) \subseteq \Kcal(X) \;\; \mbox{for all $X \subseteq \On$}. \]
In this case $\Kcal$ is called an \emph{extension} of $\Hcal$. The following observation is immediate from these definitions.

\begin{lemma} \label{l:o1}
If $\Hcal$ is a d-operator, then we have for all finite sets of ordinals $\mfrak,\nfrak$:
\begin{enumerate}[(1)]
\item
$\Hcal[\mfrak]$ is a d-operator and an extension of $\Hcal$.
\item
If $\mfrak \subseteq \Hcal(\varnothing)$, then $\Hcal[\mfrak] = \Hcal$.
\item
If $\nfrak \subseteq \Hcal[\mfrak](\varnothing)$, then $\Hcal[\nfrak] \subseteq \Hcal[\mfrak]$.
\end{enumerate}
\end{lemma}

Now we turn to specific operators $\Hcal_\sigma$. They will play a central role in the collapsing process of $\rs$ derivations; see Theorem~\ref{t:collapsing}.

\begin{definition}\rm \label{d:deriv-op}
For any $\sigma$, the operator
\[ \Hcal_\sigma \;:\; \Pow(\On) \to \Pow(\On) \]
is defined by setting for all $X \subseteq \On$:
\[ \Hcal_\sigma(X) \;:=\; 
   \bigcap\, \{C(\alpha,\beta) : X \subseteq C(\alpha,\beta) \;\,\&\;\, \sigma < \alpha \}. \]
\end{definition}

In the following lemmas we collect those properties of these operators that will be needed later. For their proofs we refer to Buchholz \cite{buchholz92a}, in particular Lemma~4.6 and Lemma~4.7.

\begin{lemma} \label{l:o2}
We have for all ordinals $\sigma,\tau$ and all $X \subseteq \On$:
\begin{enumerate}[(1)]
\item
$\Hcal_\sigma$ is a d-operator.
\item
$\Hcal_\sigma(\varnothing) = C(\sigma+1,0)$.
\item
$\tau \leq \sigma \;\;\mbox{and}\;\; \tau \in \Hcal_\sigma(X) \impl \psi(\tau) \in \Hcal_\sigma(X)$. 
\item
$\sigma < \tau \impl \Hcal_\sigma \subseteq \Hcal_\tau$.
\item
$\Hcal_\sigma(X) \cap \Omega$ is an ordinal.
\end{enumerate}
\end{lemma}

\smallskip

\begin{lemma} \label{l:o3}
Let $\mfrak$ be a finite set of ordinals and $\sigma$ be an ordinal such that the following conditions are satisfied:
\[ \mfrak \subseteq C(\sigma+1,\psi(\sigma+1)) \cap \Omega \quad\mbox{and}\quad
\sigma \in \Hcal_\sigma[\mfrak](\varnothing). \]
Given any $\alpha$ and $\beta$, we then have for $\alphah := \sigma+ \omega^{\Omega+\alpha}$ and $\betah := \sigma + \omega^{\Omega + \beta}$:
\begin{enumerate}[(1)]
\item
$\alpha \in \Hcal_\sigma[\mfrak](\varnothing) \impl \alphah \in \Hcal_\sigma[\mfrak](\varnothing)
\;\;\mbox{and}\;\; \psi(\alphah) \in \Hcal_\alphah[\mfrak](\varnothing)$.
\item
$\alpha \in \Hcal_\sigma[\mfrak](\varnothing) \;\;\mbox{and}\;\; \alpha < \beta \impl
\psi(\alphah) < \psi(\betah)$.
\item
$\Hcal_\sigma[\mfrak](\varnothing) \cap \Omega \,\subseteq\, \psi(\sigma+1)$.
\end{enumerate}
\end{lemma}

\medskip

\noindent\textbf{Convention}. From now on the letter $\Hcal$ will be used as a metavariable that ranges over d-operators.

\subsection{Operator controlled derivations in $\rs$}

For formulas $F$, finite formula sets $\Gamma$, and d-operators $\Hcal$ we define
\[ \Hcal[F] := \Hcal[|F|] \quad\mbox{and}\quad \Hcal[\Gamma] := \Hcal[|\Gamma|]. \]
Likewise, if each $\ffrak_1, \ffrak_2, \ldots, \ffrak_k$ is an ordinal, a formula or a finite set of formulas we define
\[ \Hcal[\ffrak_1,\ffrak_2] \;:=\; \Hcal[\ffrak_1][\ffrak_2], \quad
   \Hcal[\ffrak_1,\ffrak_2,\ffrak_3] \;:=\; \Hcal[\ffrak_1,\ffrak_2][\ffrak_3], \; \ldots  \]
and you can convince yourself immediately that the order of the expressions $\ffrak_1, \ffrak_2, \ldots, \ffrak_k$ does not matter.

\begin{definition}\rm \label{d:derivation}
Given a derivation operator $\Hcal$ and a finite set $\Theta$ of $\rs$ formulas, $\Hcal \prov{\sigma}{\rho}\, \Theta$ is defined by recursion on $\sigma$:
\begin{enumerate}
\item
If $\Theta$ is an axiom of $\rs$ and $\{\sigma\} \cup |\Theta| \subseteq \Hcal(\varnothing)$, then $\Hcal \prov{\sigma}{\rho}\, \Theta$.
\item
For the inference rules of $\rs$ we proceed according to the ordinal assignments listed below; in addition, we always demand that $\{\sigma\} \cup |\Theta| \subseteq \Hcal(\varnothing)$ for the conclusions $\Theta$ of all rules.
\end{enumerate}
\end{definition}

\bigskip

\noindent\textbf{Inference rules of $\rs$ with ordinal assignments}

\allowdisplaybreaks

\begin{eqnarray*}
(\lor) &
\proofrule{\Hcal \prov{\sigma_0}{\rho}\, \Gamma,\,F,\, G}
{\Hcal \prov{\sigma}{\rho}\, \Gamma,\, F \lor G}
&
\begin{array}{r}
\sigma_0 < \sigma
\end{array}
\\[4.5ex]
(\land) &
\proofrule{\Hcal \prov{\sigma_0}{\rho}\, \Gamma,\,F \qquad 
\Hcal \prov{\sigma_1}{\rho}\, \Gamma,\, G}
{\Hcal \prov{\sigma}{\rho}\, \Gamma,\, F \land G}
&
\begin{array}{r}
\sigma_0, \sigma_1 < \sigma
\end{array}
\\[4.5ex]
(\neg M) &
\proofrule{\Hcal[\beta] \prov{\sigma_\beta}{\rho}\,\Gamma,\, \neg \Mb(a) \,\,\, 
\mbox{for all $\beta < \lambda$}}
{\Hcal \prov{\sigma}{\rho} \Gamma,\, \neg M_\lambda(a)}
&
\begin{array}{r}
\sigma_\beta < \sigma
\\[0.5ex]
\mbox{$\lambda$ limit}
\end{array}
\\[4.5ex]
(\exists) &
\proofrule{\Hcal \prov{\sigma_0}{\rho}\, \Gamma,\,\Mb(a) \land F[a]}
{\Hcal \prov{\sigma}{\rho}\, \Gamma,\, \exists x F[x]}
&
\begin{array}{r}
\sigma_0 < \sigma
\\[0.5ex]
\beta < \sigma
\end{array}
\\[4.5ex]
(\forall) &
\proofrule{\Hcal[\beta] \prov{\sigma_\beta}{\rho}\, \Gamma,\, \neg \Mb(a) \lor F[a]
\,\,\,\mbox{for all $\beta < \Omega$ and all $a$}}
{\Hcal \prov{\sigma}{\rho}\, \Gamma,\, \forall x F[x]}
&
\begin{array}{r}
\beta < \sigma_\beta < \sigma
\end{array}
\\[4.5ex]
(\exists^\alpha) &
\proofrule{\Hcal \prov{\sigma_0}{\rho}\, \Gamma,\, \Mb(a) \land F[a]}
{\Hcal \prov{\sigma}{\rho}\, \Gamma,\, \exists x^\alpha F[x]}
&
\begin{array}{r}
\sigma_0 < \sigma
\\[0.5ex]
\beta \leq \alpha,\, \beta < \sigma
\end{array}
\\[4.5ex]
(\forall^\alpha) &
\proofrule{\Hcal[\beta] \prov{\sigma_\beta}{\rho}\, \Gamma,\, \neg\Mb(a) \lor F[a]\,\,\, 
\mbox{for all $\beta \leq \alpha$ and all $a$}}
{\Hcal \prov{\sigma}{\rho}\, \Gamma,\, \forall x^\alpha F[x]}
&
\begin{array}{r}
\beta < \sigma_\beta < \sigma
\end{array}
\\[4.5ex]
(b\exists) &
\proofrule{\Hcal \prov{\sigma_0}{\rho}\, \Gamma,\, b \in a \land F[b]}
{\Hcal \prov{\sigma}{\rho}\, \Gamma,\, (\exists x \in a)F[x]}
&
\begin{array}{r}
\sigma_0 < \sigma
\end{array}
\\[4.5ex]
(b\forall) &
\proofrule{\Hcal \prov{\sigma_0}{\rho}\, \Gamma,\, b \notin a \lor F[b]\,\,\, \mbox{for all $b$}}
{\Hcal \prov{\sigma}{\rho}\, \Gamma,\, (\forall x \in a)F[x]}
&
\begin{array}{r}
\sigma_0 < \sigma
\end{array}
\\[4.5ex]
(\exists_2) &
\proofrule{\Hcal \prov{\sigma_0}{\rho}\, \Gamma,\, F[U]}
{\Hcal \prov{\sigma}{\rho}\, \Gamma,\, \exists X F[X]}
&
\begin{array}{r}
\sigma_0 < \sigma
\end{array}
\\[4.5ex]
(\forall_2) &
\proofrule{\Hcal \prov{\sigma_0}{\rho}\, \Gamma,\, F[U]\,\,\, \mbox{for all $U$}}
{\Hcal \prov{\sigma}{\rho}\, \Gamma,\, \forall X F[X]}
&
\begin{array}{r}
\sigma_0 < \sigma
\end{array}
\\[4.5ex]
\cut &
\proofrule{\Hcal \prov{\sigma_0}{\rho}\, \Gamma,\,F \qquad \Hcal \prov{\sigma_0}{\rho}\, \Gamma,\, \neg F}{\Hcal \prov{\sigma}{\rho}\, \Gamma}
&
\begin{array}{r}
\sigma_0 < \sigma
\\[0.5ex]
\rk(F) < \rho
\end{array}
\\[4.5ex]
\sref &
\proofrule{\Hcal \prov{\sigma_0}{\rho}\, \Gamma,\, F}
{\Hcal \prov{\sigma}{\rho}\, \Gamma,\, \exists z F^{(z)}}
&
\begin{array}{r}
\sigma_0, \Omega < \sigma
\\[0.5ex]
F \in \Scal_0
\end{array}
\\[4.5ex]
\BC &
\proofrule{\Hcal \prov{\sigma_0}{\rho}\, \Gamma,\, F^\beta}
{\Hcal \prov{\sigma}{\rho}\, \Gamma,\, \exists z^{\beta{+}\omega}F^{(z)}}
&
\begin{array}{r}
\sigma_0  < \sigma
\\[0.5ex]
F \in \Scal_0
\end{array}
\end{eqnarray*}

\noindent This concludes the list of the rules of inference of $\rs$ with their respective ordinal assignments.

The following lemmas collect a few basic properties of the system $\rs$. They are proved by straightforward induction on $\alpha$.

\begin{lemma}[Weakening]
\[ \Hcal \prov{\sigma}{\rho_0}\, \Gamma \und \sigma \leq \tau \in \Hcal(\varnothing) \und \rho_0 \leq \rho_1 \und |\Theta| \subseteq \Hcal(\varnothing) \impl \Hcal \prov{\tau}{\rho_1}\, \Gamma, \Theta. \]
\end{lemma}

\smallskip

\begin{lemma}[Inversion] \label{l:inversion}\quad
\begin{enumerate}
\item
$\Hcal \prov{\sigma}{\rho}\, \Gamma,\, F_1 \lor F_2 \und \Omega  \leq \rk(F_1 \lor F_2) \impl 
\Hcal \prov{\sigma}{\rho}\, \Gamma,\, F_1,\, F_2$.
\item
$\Hcal \prov{\sigma}{\rho}\, \Gamma,\, F_1 \land F_2 \und \Omega  \leq \rk(F_1 \land F_2)
\und i \in \{1,2\} \impl \Hcal \prov{\sigma}{\rho}\, \Gamma,\, F_i$.
\item
$\Hcal \prov{\sigma}{\rho}\, \Gamma, \forall x F[x] \und \beta \in \Hcal(\varnothing) \cap \Omega
\impl \Hcal \prov{\sigma}{\rho}\, \Gamma,\, \neg \Mb(a) \llor F[a]$.
\item
$\Hcal \prov{\sigma}{\rho}\, \Gamma, \forall x F[x] \und \beta \in \Hcal(\varnothing) \cap \Omega
\impl \Hcal \prov{\sigma}{\rho}\, \Gamma,\, \forall x^\beta F[x]$.
\item
$\Hcal \prov{\sigma}{\rho}\, \Gamma, \forall x^\alpha F[x] \und \beta \in \Hcal(\varnothing) \und
\beta \leq \alpha
\impl \Hcal \prov{\sigma}{\rho}\, \Gamma,\, \neg \Mb(a) \llor F[a]$.
\item
$\Hcal \prov{\sigma}{\rho}\, \Gamma,\, (\forall x \in a)F[x] \und \Omega \leq \rk(F[u)])
\impl \Hcal \prov{\sigma}{\rho}\, \Gamma,\, b \notin a \llor F[s]$.
\item
$\Hcal \prov{\sigma}{\rho}\, \Gamma,\, \forall X F[X] \und \Omega \leq \rk(F[U]) \impl 
\Hcal \prov{\sigma}{\rho}\, \Gamma,\, F[U]$.
\end{enumerate}
\end{lemma}

\medskip

\begin{lemma}[Boundedness] \label{l:boundedness}
Assume $F \in \Scal$, $\sigma \leq \tau < \Omega$, and $\tau \in \Hcal(\varnothing)$. Then we have for all $\Gamma$ that
\[ \Hcal \prov{\sigma}{\rho}\, \Gamma,\, F \impl
   \Hcal  \prov{\sigma}{\rho}\, \Gamma,\, F^\tau. \]
\end{lemma}

This boundedness lemma is one of the main ingredients of impredicative cut elimination and will play a central role in the proof of the collapsing theorem in Subsection~\ref{ss:collapsing}. However, before dealing with cut elimination we turn to some embedding results.

\subsection{Embedding}

In a next step we show that $\kp + \poocas$ can be embedded into $\rs$, and as it will turn out, that finite derivations in $\kp + \poocas$ translate into uniform infinitary derivations in $\rs$. We begin with showing (TnD).

\begin{lemma}\label{l:ausvoll}
For any $\Hcal$ and any $F$ we have that
\[ \Hcal[F] \prov{\rk(F)  \natsum \rk(F)}{0}\, \neg F,\, F \]
\end{lemma}

\begin{proof}
We proceed by induction on  $\rk(F)$. If $F$ is from $\Bcal$, then $\neg F,\, F$ is an axiom, and we are done. Now suppose that $F$ is of the form $\exists x G[x]$ and observe that, for any $\beta < \Omega$, the rank of the formula $\Mb(a) \land G[a]$ is independent of $a$. Therefore, we can set
\[ \rho_\beta \;:=\; \rk(\Mb(a) \land G[a]) \quad\mbox{and}\quad \rho \;:=\; \rk(\exists x G[x]) \]
and immediately obtain $\beta \,<\, \rho_\beta \natsum \rho \,<\, \rho \natsum \rho$. Moreover,
\[ \Hcal[\Mb(a) \land G[a]] \,=\, \Hcal[\exists x G[x],\beta]. \]
Thus the induction hypothesis yields
\[ \Hcal[\exists x G[x],\beta] \prov{\rho_\beta \natsum \rho_\beta}{0}\,
\neg\Mb(a) \lor \neg G[a],\, \Mb(a) \land G[a] \]
and an application of $(\exists)$ gives us
\[ \Hcal[\exists x G[x],\beta] \prov{\rho_\beta \natsum \rho}{0}\,
\neg\Mb(a) \lor \neg G[a],\, \exists x G[x]]. \]
This is so for all $\beta < \Omega$ and all $a$. Therefore, by means of $(\forall)$ we obtain
\[ \Hcal[\exists x G[x]] \prov{\rho \natsum \rho}{0}\,
\forall x \neg G[x],\, \exists x G[x]], \]
as desired. The other cases are similar.
\end{proof}

\smallskip

\begin{lemma}[Lifting] \label{l:lifting}
For $\rho := \rk(\exists x^\alpha F[x])$  and all $\Hcal$ we have that
\[ \Hcal[\exists x^\alpha F[x]] \prov{\rho \natsum \rho}{0}\, \neg \exists x^\alpha F[x],\, \exists x F[x]. \]
\end{lemma}

\begin{proof}
For $\rho_\beta := \rk(\Mb(a) \land F[a])$ we obtain from the previous lemma that
\[ \Hcal[F[u],\beta] \prov{\rho_\beta \natsum \rho_\beta}{0}\, \neg\Mb(a) \llor \neg F[a],\, \Mb(a) \lland F[a] \]
for all $\beta \leq \alpha$. By applying  $(\exists)$ and $(\forall^\alpha)$ our assertion follows.
\end{proof}

\smallskip

\begin{lemma}[$\in$-induction]
For every formula $F[u]$ and all $\Hcal$ we have
\[ \Hcal[F[u]] \prov{\sigma+\Omega+1}{\Omega}\, \forall x((\forall y \in x)F[y] \to F[x]) \tto
\forall x F[x], \]
where $\sigma := \omega^{\rk(\forall x((\forall y \in x)F[y] \to F[x]))}$.
\end{lemma}

\begin{proof}
We set 
\[ G \;:=\; \forall x((\forall y \in x)F[y] \tto F[x]), \quad \sigma \;:=\; \omega^{\rk(G)}, \quad
   \gamma_\alpha \;:=\; \sigma \natsum \omega^\alpha \]
and prove in a first step
\[ \Hcal[F,\alpha] \prov{\gamma_\alpha}{\Omega}\, \neg G, \neg\Ma(a),\, F[a] \tag{*} \]
by induction on $\alpha$. If $|\alpha| = 0$, then simply apply Axiom~(\ref{a:5}). If $\alpha$ is a limit, then the assertion is immediate from the induction hypothesis by means of inference rule
$(\neg M)$. Now suppose $\alpha = \beta + 1$. Then $\Hcal[F,\beta] = \Hcal[F,\alpha]$ and the induction hypothesis gives us
\[ \Hcal[F,\alpha] \prov{\gamma_\beta}{\Omega}\, \neg G,\, \neg\Mb(b),\, F[b]. \]
In view of Axiom~(\ref{a:12}) we also have
\[ \Hcal[F,\alpha] \prov{0}{0}\, \neg \Ma(a),\, b \notin a,\, \Mb(b). \]
Therefore, a cut yields
\[ \Hcal[F,\alpha] \prov{\gamma_\beta +1}{\Omega}\, \neg G,\, \neg\Ma(a),\, b \notin a,\, F[b]. \]
Using $(\lor)$ and $(b\forall)$ we can continue with
\[ \Hcal[F,\alpha] \prov{\gamma_\beta +3}{\Omega}\, \neg G,\, \neg\Ma(a),\, (\forall x \in a)F[x]. \]
On the other hand, the previous lemma also implies that
\[ \Hcal[F,\alpha] \prov{\sigma}{\Omega}\, \neg F[a],\, F[a]. \]
From the previous two assertions we obtain
\[ \Hcal[F,\alpha] \prov{\gamma_\beta +4}{\Omega}\, \neg G,\, \neg\Ma(a),\, 
(\forall x \in a)F[x] \lland \neg F[a],\, F[a] \]
via $(\land)$ and from the latter
\[ \Hcal[F,\alpha] \prov{\gamma_\beta +5}{\Omega}\, \neg G,\, \neg\Ma(a),\, \Ma(a) \lland
((\forall x \in a)F[x] \lland \neg F[a]),\, F[a] \]
via Axiom~(\ref{a:1}) and $(\land)$.  Thus $(\exists)$ and the definition of $G$ lead to
\[ \Hcal[F,\alpha] \prov{\gamma_\alpha}{\Omega}\, \neg G,\, \neg\Ma(a),\, F[a]. \]
Therefore, (*) is proved. 

The rest is simple. The rule $(\lor)$ applied to (*) yields
\[ \Hcal[F,\alpha] \prov{\gamma_\alpha}{\Omega}\, \neg G,\, \neg\Ma(a) \llor F[a] \]
for all $\alpha < \Omega$ and all $a$. Hence we are in the position to apply $(\forall)$ and deduce
\[ \Hcal[F,] \prov{\sigma+\Omega}{\Omega}\, \neg G,\, \forall x F[x]. \]
From this our assertion follows by applying $(\lor)$.
\end{proof}

\smallskip

\begin{lemma}
We have for all $\Hcal$, all $a,b$, and all $\alpha,\beta$:
\begin{enumerate}[(1)]
\item
$\Hcal \prov{3}{0}\, M_1(\varnothingunder) \lland (\forall x \in \varnothingunder)(x \neq x)$.
\item
$\Hcal \prov{0}{0} M_{\omega+1}(\omegaunder)$.
\item
$\Hcal \prov{0}{0}\, \neg\Ma(a),\, (a \in \omegaunder \llra \FinOrd[a])$.
\item
$\Hcal[\alpha,\beta] \prov{\omega^{\beta+2}}{\omega\beta+\omega+\omega}\, \neg\Ma(a),\, \neg\Mb(b),\,
\exists z(a \in x \lland b \in z)$\quad for $\alpha \leq \beta$.
\item
$\Hcal[\alpha] \prov{\omega^{\alpha+2}}{\omega\alpha + \omega}\, \neg\Ma(a),\, 
\exists z(\forall y \in a)(\forall x \in y)(x \in z)$.
\end{enumerate}
\end{lemma}

\begin{proof}
These five assertions are immediate consequences of the respective axioms and Lemma~\ref{l:lifting}.
\end{proof}

It $\vec{\alpha} = \alpha_1,\ldots,\alpha_k$ and $\vec{a} = a_1,\ldots,a_k$, then we write
$\neg M_{\vec{\alpha}}(\vecc{a})$ for the set
\[ \{\neg M_{\alpha_1}(a_1),\ldots,\neg M_{\alpha_k}(a_k) \}. \]
Now we turn to $\Delta_0$ separation and our form of $\Pi^1_1$ class comprehension. Again, they are basically given by the axioms. However, for ($\Delta_0$-Sep) a cut is needed.

\smallskip

\begin{lemma}[($\Delta_0$-Sep) and  $\poocas$] \label{l:sepucomp}
Let $A[u,\vecc{v}]$ and $B[u,\vec{v},W]$ be $\Delta_0$ formulas of $\ltset$ whose free set variables are from the list $u,\vec{v}$. Then we have for all $\Hcal$, all $\alpha,\vec{\beta}$, and all $a, \vec{b}$:
\begin{enumerate}[(1)]
\item
$\Hcal[\alpha,\vecc{\beta}] \prov{\rho}{\Omega}\,
\neg \Ma(a),\, \neg M_{\vec{\beta}}(\vecc{b}),\, \exists z( z = \{x \in a : A[x,\vecc{b}]\})$, \;
where $\rho := \omega^{\max(\alpha,\vecc{\beta})+2}$.
\item
$\Hcal[\vecc{\beta}] \prov{0}{0}\,
\neg M_{\vec{\beta}}(\vecc{b}),\, 
\exists Z(\forall x \in \omegaunder)(Z(x) \llra \forall Y B[x,\vec{b},Y])$.
\end{enumerate}
\end{lemma}

\begin{proof}
The first assertion is obtained from Axiom~(\ref{a:15}) and Lemma~\ref{l:lifting} by a cut. The second is a direct consequence of Axiom~(\ref{a:16}).
\end{proof}

\begin{lemma}[$\Scal_0$ reflection] \label{l:siref}
Let $A[\vecc{u}]$ be a formula from $\Scal_0$ whose free set variables are from the list 
$\vec{u}$. Then we have for all $\Hcal$, all $\vec{\alpha}$, and all $\vec{a}$ that
\[ \Hcal[\vec{\alpha}) \prov{\sigma}{0}\, \neg M_{\vec{\alpha}}(\vecc{a}),\,  A[\vecc{a}] \tto
\exists z A^{(z)}[\vecc{a}], \]
where $\sigma \,:=\, \omega^{\rk(A[\vecc{u}]) +1}$.
\end{lemma}

\begin{proof}
In view of Lemma~\ref{l:ausvoll} we have
\[ \Hcal[\vec{\alpha}] \prov{\rho}{0}\, \neg M_{\vec{\alpha}}(\vecc{a}),\,
\neg A[\vec{a}],\, A[\vec{a}] \]
for $\rho := \rk(A[\vec{a}]) \natsum \rk(A[\vec{a}])$.  Since $A[\vec{u}]$ is a proper $\Sigma$ formula, we know that $\Omega \leq \rho$. Now we apply $\sref$ and obtain
\[ \Hcal[\vec{\alpha}] \prov{\rho+1}{0}\, \neg M_{\vec{\alpha}}(\vecc{a}),\,
\neg A[\vec{a}],\, \exists z A^{(z)}[\vec{a}] \]
Thus an application of $(\lor)$ yields our assertion.
\end{proof}

\smallskip

\begin{theorem}[Embedding] \label{t:embedding}
Let $\Gamma[u_1,\ldots,u_k]$ be a finite set of $\ltset$ formulas whose free set variables are exactly those indicated and suppose that
\[ \kp + \poocas \,\vdash^h\, \Gamma[u_1,\ldots,u_k]. \]
Then there exist $m,n < \omega$ such that
\[ \Hcal[\vecc{\alpha}] \prov{\omega^{\Omega+m}}{\Omega+n}\,
\neg M_{\vec{\alpha}}(\vecc{a}),\, \Gamma[\vecc{a}] \]
for all $\Hcal$, all $\vec{\alpha} = \alpha_1,\ldots,\alpha_k$, and all $\vecc{a} = a_1,\ldots,a_k$.
\end{theorem}

\begin{proof}
We proceed by induction on the length $h$ of the derivation of $\Gamma[\vec{u}]$ in the Tait-style formalization of $\kp + \poocas$. If $\Gamma[\vecc{u}]$ is an axiom, then the assertion follows from Lemma~\ref{l:ausvoll} -- Lemma~\ref{l:siref}. Observe that $\Delta_0$ collection is a special case of  $\Scal_0$ reflection.

\smallskip

As a first example of an inference rule we treat $(\forall)$. Then $\Gamma[\vecc{u}]$ contains a formula of the form $\forall x A[x,\vecc{u}]$ and in $\kp + \poocas$ there is a shorter derivation of
\[ \Gamma'[\vecc{u}],\, A[v,\vecc{u}], \]
where $v$ is a fresh set variable not occurring in $\Gamma[\vecc{u}]$ and $\Gamma'[\vecc{u}]$ is either $\Gamma[\vecc{u}]$ or $\Gamma[\vecc{u}] \setminus \{\forall x A[x,\vecc{u}]\}$.
Thus the induction hypothesis provides us with $m_0,n_0 < \omega$ such that
\[ \Hcal[\vec{\alpha},\beta] \prov{\omega^{\Omega+m_0}}{\Omega+n_0}\,
\neg M_{\vec{\alpha}}(\vecc{a}),\, \neg \Mb(b),\, \Gamma'[\vecc{a}],\,
A[b,\vecc{a}] \]
for all $\vec{\alpha},\beta$ and all $\vec{a},b$. Now we apply $(\lor)$ and then $(\forall)$ and obtain
\[ \Hcal[\vec{\alpha}] \prov{\omega^{\Omega+m}}{\Omega+n_0}\,
\neg M_{\vec{\alpha}}(\vecc{a}),\, \Gamma'[\vecc{a}],\, \forall xA[x,\vecc{a}] \]
for $m := m_0+1$ and arbitrary $\vec{\alpha}$ and $\vec{a}$. Since $\forall xA[x,\vecc{a}]$belongs to
$\Gamma[\vecc{a}]$ this is our assertion.

\smallskip

As second example we treat the rule $(\exists)$. Then $\Gamma[\vec{u}]$ contains a formula of the form $\exists x A[x,\vecc{u}]$ and in $\kp + \poocas$ there is a shorter derivation of
\[ \Gamma'[\vecc{u}],\, A[b,\vecc{u}], \]
for some set term $b$ and $\Gamma'[\vecc{u}]$ is either $\Gamma[\vecc{u}]$ or $\Gamma[\vecc{u}] \setminus \{\forall x A[b,\vecc{u}]\}$. Now we have to distinguish three cases:

\smallskip

\noindent (i) $b$ is the term $\varnothingunder$ or the term $\omegaunder$. Then the induction hypothesis supplies us with $m_0,n_0 < \omega$ such that
\[ \Hcal[\vecc{\alpha}] \prov{\omega^{\Omega+m_0}}{\Omega+n_0}\,
\neg M_{\vec{\alpha}}(\vecc{a}),\, \Gamma[\vec{a}],\, A[b,\vecc{a}] \]
for all $\vec{\alpha}$ and all $\vec{a}$. Together with Axiom~(\ref{a:6}) or Axiom~(\ref{a:8}) and $(\land)$ we thus obtain 
\[ \Hcal[\vecc{\alpha}] \prov{\omega^{\Omega+m_0}+1}{\Omega+n_0}\,
\neg M_{\vec{\alpha}}(\vecc{a}),\, \Gamma[\vec{a}],\, \Mb(b) \land A[b,\vecc{a}] \]
(where $\beta$ is $1$ or $\omega+1$) and thus, by $(\exists)$,
\[ \Hcal[\vecc{\alpha}] \prov{\omega^{\Omega+m}}{\Omega+n_0}\,
\neg M_{\vec{\alpha}}(\vecc{a}),\, \Gamma[\vec{a}],\, \exists x A[x,\vecc{a}]. \]
for $m := m_0+1$ and arbitrary $\vec{\alpha}$ and $\vec{a}$. 

\smallskip

\noindent(ii) $b$ is the variable $u_i$, $1 \leq i \leq n$. Then the induction hypothesis yields
\[ \Hcal[\vecc{\alpha}] \prov{\omega^{\Omega+m_0}}{\Omega+n_0}\,
\neg M_{\vec{\alpha}}(\vecc{a}),\, \Gamma[\vec{a}],\, A[a_i,\vecc{a}] \]
for suitable $m_0,n_0 < \omega$, all $\vec{\alpha}$ and all $\vec{a}$. Now Axiom~(\ref{a:1}) and $(\land)$ give us
\[ \Hcal[\vecc{\alpha}] \prov{\omega^{\Omega+m_0}+1}{\Omega+n_0}\,
\neg M_{\vec{\alpha}}(\vecc{a}),\, \Gamma[\vec{a}],\, M_{\alpha_i}(a_i)  \land A[a_i,\vecc{a}] \]
and from there we proceed as in (i).

\smallskip

\noindent(iii) $b$ is a set  variable not occurring in $\Gamma[\vecc{u}]$. Then the induction hypothesis supplies us with $m_0,n_0 < \omega$ such that
\[ \Hcal[\vecc{\alpha}] \prov{\omega^{\Omega+m_0}}{\Omega+n_0}\,
\neg M_{\vec{\alpha}}(\vecc{a}),\,\neg M_1(\varnothingunder),\,
\Gamma[\vec{a}],\, A[\varnothingunder,\vecc{a}] \]
for all $\vec{\alpha}$ and all $\vec{a}$. Now we make first use of Axiom~(\ref{a:6}) and a cut and then of Axiom~(\ref{a:6}) and $(\land)$ to derive 
\[ \Hcal[\vecc{\alpha}] \prov{\omega^{\Omega+m_0}+2}{\Omega+n_0}\,
\neg M_{\vec{\alpha}}(\vecc{a}),\, \Gamma[\vec{a}],\, M_1(\varnothingunder) \land A[\varnothingunder,\vecc{a}] \]
from where we can proceed again as in (i).

\smallskip

In all three cases we obtain
\[ \Hcal[\vecc{\alpha}] \prov{\omega^{\Omega+m}}{\Omega+n_0}\,
\neg M_{\vec{\alpha}}(\vecc{a}),\, \Gamma[\vec{a}],\, \exists x A[x,\vecc{a}] \]
for $m := m_0+1$ and arbitrary $\vec{\alpha}$ and $\vec{a}$. Since $\exists x A[x,\vecc{a}]$ belongs to $\Gamma[\vecc{a}]$ this finishes $(\exists)$.

\smallskip

All other cases are straightforward from the induction hypothesis or can be treated in the same vein. 
\end{proof}

\subsection{Predicative cut elimination}

The rules of inference of $\rs$ can be divided into two classes: (i) In all rules except $\sref$ the principal formula is more complex than the corresponding formula(s) in the premise(s). We may, therefore, consider these rules as \emph{predicative} rules. (ii) The rule $\sref$, on the other hand, transforms (in the general case) a formula of rank greater than $\Omega$ into one of rank $\Omega$. This is a sort of impredicativity and we consider $\sref$ as an \emph{impredicative} inference. 

Also remember that the principal formulas of the axioms of $\rs$ belong to $\Bcal$ and, therefore, have rank less than $\Omega$. These formulas are an obstacle for cut elimination below $\Omega$. However, above $\Omega$ we can follow the pattern of the usual (predicative) cut elimination.

In this section we sketch how all cuts with cut formulas of ranks greater than $\Omega$ can be eliminated by standard methods as presented, for example, in Sch\"utte \cite{schuette77} or Buchholz \cite{buchholz92a}.

\begin{lemma}
Let $F$ be a formula of the form $G_1 \lor G_2$, $(\exists x \in a)G[x]$, $\exists x^\alpha G[x]$,
$\exists x G[x]$ or $\exists X G[X]$ and assume that the rank $\rho$ of $F$ is greater than $\Omega$. Then we have:
\[ \Hcal \prov{\alpha}{\rho}\, \Gamma, \neg F \;\;\mbox{and}\;\;\; \Hcal \prov{\beta}{\rho}\, \Theta, F \impl
\Hcal \prov{\alpha + \beta}{\rho} \Gamma,\, \Theta. \]
\end{lemma}

\begin{proof}
By straightforward induction on $\beta$.
\end{proof}

\bigskip

\begin{theorem}[Predicative cut elimination] \label{t:pred-cut-el}
We have the following two elimination results, where the second is an immediate consequence of the first. For all $\Hcal$, $\Gamma$, $\alpha$, and all $n < \omega$:
\begin{enumerate}[(1)]
\item
$\Hcal \prov{\alpha}{\Omega + n +2}\, \Gamma \impl 
\Hcal \prov{\omega^\alpha}{\Omega + n + 1}\, \Gamma$.
\item
$\Hcal \prov{\alpha}{\Omega + n +1}\, \Gamma \impl 
\Hcal \prov{\omega_n[\alpha]}{\Omega + 1}\, \Gamma$.
\end{enumerate}
Recall that $\omega_0[\alpha] := \alpha$ and  $\omega_{n+1}[\alpha] := \omega^{\omega_n[\alpha]}$.
\end{theorem}

\begin{proof}
The proof of the first assertion is standard. For details see, for example, Buchholz \cite{buchholz92a}.  The second assertion is an immediate consequence of the first. 
\end{proof}

\subsection{Collapsing} \label{ss:collapsing}

We begin  this subsection by showing that specific operator controlled derivations of finite sets of formulas from $\Scal \cup \Bcal$ in which all cut formulas have ranks $\leq \Omega$ can be collapsed into derivations of depth and cut rank less than $\Omega$. This technique -- called collapsing technique -- is a corner stone of \emph{impredicative proof theory}. See Buchholz \cite{buchholz92a} for more on general collapsings.

Together with the results of the previous subsection this collapsing theorem will then lead to our main theorem about the $\Scal \cup \Bcal$ fragment of $\rs$ derivations and the $\Sigma$ formulas provable in $\kp + \poocas$.

For the collapsing theorem we work with the derivation operators $\Hcal_\sigma$ introduced in Definition~\ref{d:deriv-op}.

\bigskip

\begin{theorem}[Collapsing] \label{t:collapsing}
Let $\Gamma$ be a finite set of formulas from $\Scal \cup \Bcal$ and suppose that 
\[ |\Gamma| \subseteq C(\sigma+1,\psi(\sigma+1)) \quad\mbox{and}\quad
\sigma \in \Hcal_\sigma[\Gamma](\varnothing). \]
Then we have, for all $\alpha$,
\[ \Hcal_\sigma[\Gamma] \provv{\alpha}{\Omega+1}\, \Gamma  \impl
\Hcal_\alphah[\Gamma] \provv{\psi(\alphah)}{\psi(\alphah)}\, \Gamma, \]
where $\alphah := \sigma + \omega^{\Omega+\alpha}$.
\end{theorem}

\begin{proof}\setcounter{equation}{0}
By induction on $\alpha$. Assume $\Hcal_\sigma[\Gamma] \prov{\alpha}{\Omega+1}\, \Gamma$. Then we certainly have (see Lemma~\ref{l:o2} and Lemma~\ref{l:o3}):
\begin{gather}
   |\Gamma | \cup \{\alpha\} \;\subseteq\; \Hcal_\sigma[\Gamma](\varnothing) \;\subseteq\; 
   \Hcal_\alphah[\Gamma](\varnothing),
   \\[1.5ex]
   \alphah \in \Hcal_\sigma[\Gamma](\varnothing) \quad\mbox{and}\quad
   \psi(\alphah) \in \Hcal_\alphah[\Gamma](\varnothing).
\end{gather}
Now we distinguish cases according to the last inference of $\Hcal_\sigma[\Gamma] \provv{\alpha}{\Omega+1}\, \Gamma$ and note that this cannot be $(\forall)$. In the following we confine our attention to the interesting cases; all others can be dealt with in a similar manner.

\medskip

\noindent (i) $\Gamma$ is an axiom. Then 
$\Hcal_{\hat{\alpha}}[\Gamma] \provv{\psi(\hat{\alpha})}{\psi(\hat{\alpha})}\, \Gamma$ follows from (1) and (2).

\medskip

\noindent (ii) The last inference was $(\forall^\tau)$. Then $\Gamma$ contains a formula of the form $\forall x^\tau F[x]$ and we have
\begin{gather}
   \Hcal_\sigma[\Gamma,\xi] \provv{\alpha_\xi}{\Omega+1}\, \Gamma,\, 
   \neg M_\xi(a) \llor F[a]
\end{gather}
with $\alpha_\xi < \alpha$ for all $\xi \leq \tau$ and all $a$. We know $\tau  \in \Hcal_\sigma[\Gamma](\varnothing)$. Lemma~\ref{l:o1}, Lemma~\ref{l:o2}(5), and (1) yield
\begin{gather}
   \xi \in \Hcal_\sigma[\Gamma](\varnothing) \subseteq \Hcal_\alphah[\Gamma](\varnothing),
   \quad
   \Hcal_\sigma[\Gamma,\xi] = \Hcal_\sigma[\Gamma],
   \quad
   \Hcal_\alphah[\Gamma,\xi] = \Hcal_\alphah[\Gamma].
\end{gather}
Hence (3) gives us
\[ \Hcal_\sigma[\Gamma] \provv{\alpha_\xi}{\Omega+1}\, \Gamma,\, 
   \neg M_\xi(a) \llor F[a] \quad\mbox{and}\quad \alpha_\xi \in \Hcal_\sigma[\Gamma]. \]
By the induction hypothesis we obtain that
\[ \Hcal_\alphaxih[\Gamma] \provv{\psi(\alphaxih)}{\psi(\alphaxih)}\, \Gamma,\, 
   \neg M_\xi(a) \llor F[a] \]
and, therefore,
\begin{gather}
   \Hcal_\alphah[\Gamma] \provv{\psi(\alphaxih)}{\psi(\alphaxih)}\, \Gamma,\, 
   \neg M_\xi(a) \llor F[a],
\end{gather}
always for all $\xi \leq \tau$ and all $a$. For $\xi \leq \tau$ we have  $\alpha_\xi < \alpha$ and $\alpha_\xi \in \Hcal_\sigma[\Gamma]$, hence Lemma~\ref{l:o3} yields $\psi(\alpha_\xi) < \psi(\alpha)$. Furthermore, since $\tau \in |\Gamma| \subseteq C(\sigma+1,\psi(\sigma+1))$ we also have
\[ \xi \leq \tau < \psi(\sigma+1) \leq \psi(\alphaxih). \]
Then $\Hcal_\alphah[\Gamma] \provv{\psi(\alphah)}{\psi(\alphah)}\, \Gamma$ follows immediately by an application of $(\forall^\tau)$.

\medskip

\noindent (iii) The last inference was $(\exists)$. Then $\Gamma$ contains a formula of the form $\exists x F[x]$ and we have
\begin{gather}
   \Hcal_\sigma[\Gamma] \provv{\alpha_0}{\Omega+1}\, \Gamma,\, \Mb(a) \lland F[a]
\end{gather}
for some $\alpha_0 < \alpha$ and some $\beta < \alpha$. In view of (6) we have
$\alpha_0,\beta \in \Hcal_\sigma[\Gamma](\varnothing)$. The induction hypothesis yields
\[ \Hcal_\alphazh[\Gamma] \provv{\psi(\alphazh)}{\psi(\alphazh)}\, \Gamma,\,
\Mb(a) \lland F[a], \]
thus also
\[ \Hcal_\alphah[\Gamma] \provv{\psi(\alphazh)}{\psi(\alphazh)}\, \Gamma,\,
\Mb(a) \lland F[a]. \]
Now recall from Lemma~\ref{l:o3}(3) that
\[ \Hcal_\sigma[\Gamma](\varnothing) \cap \Omega \,\subseteq\, \psi(\sigma+1) \,\leq\,
\psi(\alphah) \]
and conclude that $\beta < \psi(\alphah)$. In addition, $\psi(\alphazh) < \psi(\alphah)$ follows from Lemma~\ref{l:o3}(2). Therefore, an application of $(\exists)$ yields
$\Hcal_\sigma[\Gamma] \provv{\psi(\alphah)}{\psi(\alphah)}\, \Gamma$.

\medskip

\noindent (iv) The last inference was $\sref$. Then there exist an $F \in \Scal_0$ and an $\alpha_0 < \alpha$ such that $\exists z F^{(z)}$ belongs to $\Gamma$ and
\begin{gather}
   \Hcal_\sigma[\Gamma] \provv{\alpha_0}{\Omega+1}\, \Gamma,\, F.
\end{gather}
The induction hypothesis immediately yields
\[ \Hcal_\alphazh[\Gamma] \prov{\psi(\alphazh)}{\psi(\alphazh)}\, \Gamma,\, F.  \]  
Since $F \in \Scal_0$, we are in the position to make use of Lemma~\ref{l:boundedness} and obtain
\begin{gather*}
   \Hcal_\alphazh[\Gamma] \provv{\psi(\alphazh)}{\psi(\alphazh)}\, \Gamma,\, F^{\psi(\alphazh)}.
\end{gather*}
Now we apply $\BC$ to derive
\begin{gather}
   \Hcal_\alphazh[\Gamma] \provv{\psi(\alphazh)+1}{\psi(\alphazh)}\, \Gamma,\,
   \exists z^\beta F^{(z)}
\end{gather}
for $\beta := \psi(\alphazh) + \omega$. On the other hand, Lemma~\ref{l:lifting} yields
\begin{gather}
\Hcal_\alphazh[\Gamma] \provv{\rho \natsum \rho}{0}\,
\neg \exists z^\beta F^{(z)},\, \exists z F^{(z)}
\end{gather}
for $\rho := \rk(\exists z^\beta F^{(z)})$. Furthermore, since $F$ is from $\Scal_0$ we can easily see that the rank of $\exists z^\beta F^{(z)}$ is smaller than $\omega\beta+\omega$.
Hence some trivial ordinal computations yield that 
$\psi(\alphazh)+1, \rho \natsum \rho < \psi(\alphah)$. Therefore, a cut applied to (8) and (9) implies $\Hcal_\alphah[\Gamma] \provv{\psi(\alphah)}{\psi(\alphah)}\, \Gamma$.

\medskip

\noindent (v) The last inference was $\cut$. Then there exist a formula $F$ and an ordinal $\alpha_0 < \alpha$ such that $\rk(F) < \Omega+1$ and
\begin{gather}
   \Hcal_\sigma[\Gamma] \provv{\alpha_0}{\Omega+1}\, \Gamma,\, F \quad\mbox{and}\quad
   \Hcal_\sigma[\Gamma] \provv{\alpha_0}{\Omega+1}\, \Gamma,\, \neg F.
\end{gather}
Therefore, $\alpha_0 \in \Hcal_\sigma[\Gamma](\varnothing)$, $|F| \subseteq \Hcal_\sigma[\Gamma](\varnothing) \cap \Omega$ and, by Lemma~\ref{l:o3}(3),
\[ |F| \,\subseteq\, \Hcal_\sigma[\Gamma](\varnothing) \cap \Omega \,\subseteq\, \psi(\sigma+1) \,\leq\, \psi(\alphah). \]
We distinguish two cases.

\medskip

\noindent (v.1) $\rk(F) < \Omega$. Then $\rk(F) < \psi(\alphah)$ according to the line above. Furthermore, $\Gamma \cup \{F,\neg F\} \subseteq \Scal \cup \Bcal$ and thus the induction hypothesis yields
\begin{gather*}
   \Hcal_\alphazh[\Gamma] \provv{\psi(\alphazh)}{\psi(\alphazh)}\, \Gamma,\, F \quad\mbox{and}\quad
   \Hcal_\alphazh[\Gamma] \provv{\psi(\alphazh)}{\psi(\alphazh)}\, \Gamma,\, \neg F.
\end{gather*}
As above, we can easily convince ourselves that $\psi(\alphazh) < \psi(\alphah)$. Hence, $\cut$ yields
$\Hcal_\alphah[\Gamma] \provv{\psi(\alphah)}{\psi(\alphah)}\, \Gamma$.

\medskip

\noindent (v.2) $\rk(F) = \Omega$. Then $F$ is of the form $\exists xG[x]$ or $\forall x G[x]$ and 
$G[u]$ is a $\Delta_0$ formula of $\ltset$. We assume that $F$ is $\exists x G[x]$. Clearly, $\exists x G[x]$ belongs to $\Scal_0$. Hence, the induction hypothesis applied to the left hand side of (10) yields
\begin{gather*}
   \Hcal_\alphazh[\Gamma] \provv{\beta}{\beta}\, \Gamma,\, \exists x G[x]
\end{gather*}
for $\beta := \psi(\alphazh)$. Since $\beta \in \Hcal_\alphazh[\Gamma](\varnothing)$ (see Lemma~\ref{l:o2}), we can apply Lemma~\ref{l:boundedness} and obtain
\begin{gather}
   \Hcal_\alphazh[\Gamma] \provv{\beta}{\beta}\, \Gamma,\, \exists x^\beta G[x].
\end{gather}
Furthermore, applying Lemma~\ref{l:inversion} to the right hand side of (10) yields
\begin{gather}
   \Hcal_\alphazh[\Gamma] \provv{\alpha_0}{\Omega+1}\, \Gamma,\, \forall x^\beta \neg G[x].
\end{gather}
Now we exploit the fact that $\forall x^\beta \neg G[x]$ is from $\Bcal$ and convince ourselves that the induction hypothesis can be applied to (12). So we obtain
\begin{gather*}
   \Hcal_\gamma[\Gamma] \provv{\psi(\gamma)}{\psi(\gamma)}\, \Gamma,\, \forall x^\beta \neg G[x]
\end{gather*}
with $\gamma = \alphazh + \omega^{\Omega+\alpha_0} = 
\sigma + \omega^{\Omega+\alpha_0} + \omega^{\Omega+\alpha_0} < \sigma + \omega^{\Omega + \alpha} = \alphah$. Moreover, it is easy to check that $\gamma \in \Hcal_\sigma[\Gamma](\varnothing)$, hence $\psi(\gamma) < \psi(\alphah)$ according to Lemma~\ref{l:o3}. An easy computation, similar to above, also shows that
\[\rk(\exists x^\beta G[x]) < \psi(\alphah). \]
Therefore, $\cut$ applied to (11) and (12) establishes
$\Hcal_\alphah[\Gamma] \provv{\psi(\alphah)}{\psi(\alphah)}\, \Gamma$.

\medskip

\noindent (vi) All other cases are trivial or can be dealt with similarly. This completes the proof of collapsing.
\end{proof}

\subsection{Finite bounds for the lengths of formulas in infinitary derivations} \label{ss:finite-bounds}

Clearly, every proof $\Pcal$ in the theory $ \kp + \poocas$ is finite and, therefore, there exists a natural number $p$ such that every formula in $\Pcal$ has a length less than $p$. We sketch now how this bound $p$ carries over to the ``proof-theoretically treated'' derivation stemming from $\Pcal$. First we introduce a suitable definition of length of an $\ls$ formula.

\begin{definition}\rm
The \emph{length} $\ell(F)$ of an $\ls$ formula $F$is inductively defined as follows:
\begin{enumerate}
\item 
If $F$ is of the form $(a \in b)$, $(a \notin b)$, $U(a)$, $\neg U(a)$, $\Ma(a)$ or $\neg \Ma(a)$, then $\ell(F) := 0$.
\item
If $F$ is of the form $(G \lor H)$ or $(G \land H)$, then $\ell(F) := \max(\ell(G),\ell(H)) +1$.
\item 
If $F$ is of the form $\exists x G[x]$, $\forall x G[x]$, $\exists x^\alpha G[x]$, $\forall x^\alpha G[x]$, $(\exists x  \in a)G[x]$ or $(\forall x  \in a)G[x]$, then $\ell(F) := \ell(G[u]) +1$.
\item
If $F$ is of the form $\exists XG[X]$ or $\forall XG[X]$, then $\ell(F) := \ell(G[U]) +1$.
\end{enumerate}
\end{definition}

Based on this definition we now introduce a refined notion of derivability in $\rs$. Given a natural number $p$ we let 
\[  \Hcal \provv{\alpha }{\rho,\,p}\, \Gamma \]
be defined as $\Hcal \provv{\alpha }{\rho}\, \Gamma$ in Definition~\ref{d:derivation}  with the additional requirements that $\ell(F) < p$ for every element $F$ of $\Gamma$.

The proofs of the following three theorems do not raise any questions of principle. We simply follow the original proofs and check by a tedious case to case analysis that the additional requirements are satisfied. In the infinitary derivations new formulas may occur; however, their lengths are always bounded by the lengths of formulas in the finite derivation in $\kp + \poocas$.

\begin{theorem}[Refined embedding] \label{t:refined-embedding}
Let $\Gamma[u_1,\ldots,u_k]$ be a finite set of $\ltset$ formulas whose free set variables are exactly those indicated and suppose that 
\[ \kp + \poocas \,\vdash^h\,\Gamma[u_1,\ldots,u_k] \]
such that $\ell(F) < p$ for all formulas occurring in this proof. Then there exist $m,n < \omega$ such that
\[ \Hcal[\vecc{\alpha}] \prov{\omega^{\Omega+m}}{\Omega+n,\,p}\,
\neg M_{\vec{\alpha}}(\vecc{a}),\, \Gamma[\vecc{a}] \]
for all $\Hcal$, all $\vec{\alpha} = \alpha_1,\ldots,\alpha_k$, and all $\vecc{a} = a_1,\ldots,a_k$.
\end{theorem}

\begin{theorem}[Refined predicative cut elimination] \label{t:refined-pred-cut-el}
For all $\Hcal$, $\Gamma$, $\alpha$, and all $n,p < \omega$,
\[ 
\Hcal \prov{\alpha}{\Omega + n +1,\,p}\, \Gamma \impl 
\Hcal \prov{\omega_n(\alpha)}{\Omega + 1,\,p}\, \Gamma. \]
\end{theorem}

\begin{theorem}[Refined collapsing] \label{t:refined-collapsing}
Let $\Gamma$ be a finite set of formulas from $\Scal \cup \Bcal$ and suppose that 
\[ |\Gamma| \subseteq C(\sigma+1,\psi(\sigma+1)) \quad\mbox{and}\quad
\sigma \in \Hcal_\sigma[\Gamma](\varnothing). \]
Then we have, for all $\alpha$ and all $p < \omega$,
\[ \Hcal_\sigma[\Gamma] \provv{\alpha}{\Omega+1,\,p}\, \Gamma  \impl
\Hcal_\alphah[\Gamma] \provv{\psi(\alphah)}{\psi(\alphah),\,p}\, \Gamma, \]
where $\alphah := \sigma + \omega^{\Omega+\alpha}$.
\end{theorem}

\section{Reduction to $\TT$} \label{s:reduction}

The goal of this section is to interpret certain $\rs$-derivations of length and cut rank 
$< \Omega$ in the theory $\TT$. Here $\pooca$ stands for the usual system of second order arithmetic with $\Pi^1_1$ comprehension and full induction on the natural numbers; it is formulated in the language $\ltwo$. This is a familiar theory in proof theory and, therefore, we  refrain from saying more about this theory. If you are interested in all details you may consult, for example Simpson \cite{simpson09}.

In the sequel -- when arguing in the framework of second order arithmetic -- we assume an arithmetization of the notation system $C(\varepsilon_{\Omega+1},0)$ such that all relevant ordinal sets, functions and relations become primitive recursive. Also, we will identify the ordinal notations with their arithmetical codes. In addition, we write $\prec$ for the primitive recursive relation defined on $C(\varepsilon_{\Omega+1},0) \cap \Omega$ by
\[ \alpha \prec \beta \;:=\; \alpha, \beta \in C(\varepsilon_{\Omega+1},0) \lland
   \alpha < \beta < \Omega. \]
Recall that (the codes for) the elements of $C(\varepsilon_{\Omega+1} \cap \Omega$ denote exactly all ordinals less than the Bachmann-Howard ordinal $\BH$. In the context of $\ltwo$, the quantifiers $\exists \alpha,\forall \beta$ range over the latter set. 

By $\mathit{TI}({<}\BH)$ we mean the scheme of transfinite induction over all initial segments (indexed externally) below the Bachmann-Howard ordinal $\BH$. More precisely, 
$\mathit{TI}({<}\BH)$ is the collection of all formulas
\[ \forall \alpha((\forall \beta \prec \alpha)\Fcal[\beta] \to \Fcal[\alpha]) \tto 
(\forall \alpha \prec \psi(\omega_n[\Omega+1]))\Fcal[\alpha] \]
where $\Fcal[u]$ is any formula of $\ltwo$ and $n$ is any natural number (or rather the $n^{th}$ numeral). 

\subsection{$\alpha$-trees} \label{s:suitable}

To interpret set theory in  $\TT$ we use well-founded trees, also called {\em suitable} trees. We will mostly follow the terminology and presentation in Simpson \cite[VII.3]{simpson09}.

\begin{definition}\label{MR1} \rm
Let variables $\sigma,\tau,\rho,\sigma_0,\sigma_1,\ldots$ range over finite sequences of naturals, i.e., elements of $\Nbb^{<\Nbb}$. 
\begin{enumerate}
\item
$\langle\rangle$ will denote the empty sequence of $\Nbb^{<\Nbb}$ and $\langle n_0,\ldots, n_r\rangle$ the sequence of numbers  $n_0,\ldots, n_r$ coded
as a single number (see \cite[II.2.6]{simpson09}).
\item
The \emph{length} of a sequence is defined by  $\mathrm{lh}(\langle\rangle)=0$ and $\mathrm{lh}(\langle n_0,\ldots, n_r\rangle)=r+1$.
\item
The \emph{concatenation} operation on $\mathbb{N}^{<\mathbb N}$ will be denoted by $\sigma\star\tau$. This means that $\sigma\star\langle\rangle=\langle\rangle\star\sigma=\sigma$
and $\langle n_0,\ldots,n_r\rangle \star\langle m_0,\ldots,m_s\rangle=\langle n_0,\ldots,n_r,m_0,\ldots,m_s\rangle$. 
\item
We write $\sigma \subset \tau$ to mean that $\sigma$ is a proper \emph{initial segment} of $\tau$, i.e., $\tau=\sigma\star \rho$ for some $\rho\ne \langle\rangle$.
\item
For a function $f\colon \mathbb{N}\to \mathbb{N}$ let $f[0]=\langle \rangle$ and $f[n+1]:=\langle f(0),\ldots,f(n)\rangle$.
\end{enumerate}
\end{definition} 

\begin{definition}\label{MR2}\rm
$T \subseteq \Nbb^{<\Nbb}$ is said to be a \emph{tree} if $T$ is nonempty and closed under initial segments, i.e., 
$(\forall \sigma \in T)(\forall \tau\subset \sigma)(\tau \in T)$.  $T$ is a \emph{suitable tree} if $T$ is also well-founded, i.e.,
\[ (\forall f : \Nbb \to \Nbb)\exists n(f[n]\notin T). \]
\end{definition}

\begin{remark}\rm
If $T$ is a tree and $\sigma\in T$, we put
\[ T^{\sigma}=\{\tau\mid \sigma\star\tau \in T\} \]
and note that if $T$ is suitable tree then so is $T^{\sigma}$.
\end{remark}

Suitable trees furnish a way of talking about sets in the language of second order arithmetic. The following definition, also taken from  Simpson \cite[VII.3]{simpson09}, tells us how to define the equality relation and element relation on suitable trees.

\begin{definition}[$=^*$ and $\in^*$]\label{MR3}\rm \quad
\begin{enumerate}
\item
Let $T$ be a tree. We write $\mathit{Iso}(X,T)$ to mean that $X \subseteq T \times T$ and, for all $\sigma,\tau \in T$,  $(\sigma,\tau)\in X$ if and only if 
\[ \forall n(\sigma\star \langle n \rangle \in T \tto 
   \exists m((\sigma\star \langle n\rangle,\tau\star\langle m\rangle)\in X)) \]
   
\vspace{-3ex}and\vspace{-4ex}

\[ \forall m(\tau\star \langle m \rangle \in T \tto
   \exists n((\sigma\star \langle n\rangle,\tau\star\langle m\rangle)\in X)). \]
\item
For  trees $S$ and $T$, we define
\[ S\oplus T \;:=\; \{\langle \rangle\} \,\cup\, \{\langle 0\rangle \star \sigma\mid \sigma\in S\}
   \,\cup\, \{\langle 1\rangle \star \tau\mid \tau\in T\} \]

\vspace{-3ex}and set\vspace{-5ex}

\begin{eqnarray*}
S=^*T &\mbox{ iff }& \exists X(\mathit{Iso}(X,S\oplus T) \lland 
(\langle 0\rangle,\langle 1\rangle)\in X),
\\[1ex]
 S\in^*T &\mbox{ iff }& \exists X(\mathit{Iso}(X,S\oplus T) \lland 
 \exists n((\langle 0\rangle,\langle 1,n\rangle)\in X)).
 \end{eqnarray*}
\end{enumerate}
\end{definition}

\begin{lemma}[$\atro$]\label{MR4}  
Let $T$ be a suitable tree. Then there exists a unique set $X$ such that $\mathit{Iso}(X,T)$. Moreover,  for all $\sigma,\tau\in T$,
\begin{eqnarray*} 
T^{\sigma}=^* T^{\tau} &\mbox{ iff }& (\sigma,\tau)\in X,
\\[1ex]
T^{\sigma}\in ^* T^{\tau} &\mbox{ iff }& \exists n((\sigma,\tau\star\langle n\rangle)\in X).\end{eqnarray*}
In particular, $X$ is an equivalence relation on $T$.
\end{lemma}

\begin{proof}
See \cite[Lemma VII.3.17]{simpson09}.
\end{proof}

\begin{corollary}[$\atro$]\label{MR5} 
Given suitable trees $S$ and $T$, one has 
\begin{eqnarray*}
S\in^* T &\mbox{ iff }& \exists n(S=^*T^{\langle n\rangle}).
\end{eqnarray*}
\end{corollary}

\begin{definition}\label{MR6}\rm
In order to model the set-theoretic naturals and $\omega$ via suitable trees we introduce trees $n^*$ for $n\in\Nbb$.
Let 
\begin{eqnarray*}
0^*&=&\{\langle\rangle\}
\\[1ex]
(n+1)^* &=& n^*\,\cup\,\{\langle n\rangle \star \sigma\mid \sigma\in n^*\}
\\[1ex]
\omega^* &=&\{\langle\rangle\}\,\cup\,\{\langle n\rangle\star\sigma\mid \sigma\in n^*,\; n\in \Nbb\}.
\end{eqnarray*}
\end{definition}

Note that $n^ *\in^* \omega^*$ and $m^*\in^*n^*$ whenever $m<n$. Since the sequences in $n^*$, and thus in $\omega^*$, are strictly descending sequences of naturals it is clear that $n^*$ and $\omega^*$ are suitable trees.
Moreover, there is a map $f\colon \omega^*\to \omega{+}1$ such that
\[ \tau\subset\sigma \lland \sigma\in\omega^*  \impl f(\sigma)\prec f(\tau).\]
 To see this, let $f(\langle\rangle)=\omega$ and $f(\langle n_1,\ldots,n_r\rangle)=n_r$. This means that $\omega^*$ is an $\omega{+}1$-tree in the sense of Definition~\ref{MR7} below.

The natural translation of the set-theoretic language into that of second order arithmetic,
$\mathsf{L}_2$, proceeds by letting quantifiers range over suitable trees and interpreting $\in$ and $=$ as $\in^*$ and $=^*$, respectively. In \cite[Theorem VII.3.22]{simpson09} it is shown that this yields an interpretation of the axioms of $\atroset$ in $\atro$.

In what follows it is our goal to interpret the richer language of $\rs$ in  $\ltwo$. For the predicates  $\Ma$ and the quantifiers $\forall x^{\alpha}$ and $\exists x^{\alpha}$ we use the notion of an $\alpha$-tree.

\begin{definition}\label{MR7}\rm
For the notation $\alpha$ of an ordinal less than $\BH$ , we say that a tree $T$ is an \emph{$\alpha$-tree} if there exists a function $f : T\to \{\beta\in \BH : \beta\prec \alpha\}$ such that
\[ (\forall \sigma \in T)\forall \tau(\tau\subset \sigma \to f(\sigma)\prec  f(\tau)]. \]
\end{definition}

The interpretation of the predicate $M_{\alpha}$ for trees $T$ will be that $M_{\alpha}(T)$ means that $T$ is an $\alpha$-tree. Quantifiers $\forall x^{\alpha}(\ldots x\ldots)$ and $\exists x^{\alpha}(\ldots x\ldots)$ will 
then be interpreted as $\forall T(T\mbox{ $\alpha$-tree} \tto \ldots T\ldots)$ and
$\exists T(T\mbox{ $\alpha$-tree} \lland \ldots T\ldots)$, respectively.

Relation variables $U,VX,Y,\ldots$ will be interpreted to range over $\omega{+}1$-trees $T$ such that $T\subseteq^*\omega^*$, i.e.,
$(\forall \langle n\rangle\in T)(T^{\langle n\rangle}\in^*\omega^*)$. It is perhaps worth mentioning  that any suitable tree satisfying $T\subseteq^*\omega^*$ can be seen to be an $\omega{+}1$-tree. 

Note that if we have the notation $\alpha$ of an ordinal less than $\Omega$ (given externally), then this ordinal is less than $\BH$ and $\TT$ proves transfinite induction along the ordinal notations $\prec \alpha$, and consequently the notion of an $\alpha$-tree becomes $\Delta^1_1$ in $\TT$. Even a slightly stronger assertion can be proved.

\begin{lemma}\label{MRalpha} 
The notion of an $\alpha$-tree is $\Delta^1_1$  in the theory $\atro+\mathit{WO}(\alpha)$, where $\mathit{WO}(\alpha)$ expresses that $\prec$ restricted to $\{\beta : \beta\prec \alpha\}$ is a wellordering.
\end{lemma}

\begin{proof}
For a tree $T$, define a function $g$ with domain $\{\beta : \beta\prec \alpha\}$ by arithmetical transfinite recursion 
as follows
\[ g(\beta) \;=\; \{\sigma\in T : (\forall \tau\in T)(\sigma\subset \tau \tto \tau\in \bigcup_{\xi\prec \beta}g(\xi))\}. \]
Note that $g$ is uniquely determined by $T$. One can then show that:
\begin{eqnarray*}\label{MR1.8.8} 
\mbox{$T$ is an $\alpha$-tree } &\Leftrightarrow & \mbox{ $T$ is the image of $g$}. \end{eqnarray*}
If $T$ is an $\alpha$-tree witnessed by a function  $f : T\to \{\beta : \beta\prec \alpha\}$, then one shows, using $\prec$-induction on $f(\sigma)$,  that $\sigma\in g(f(\sigma))$ for all $\sigma\in T$, hence $T$ is the image of $g$.

 Conversely, if $T$ is the image of $g$, then a function $f : T\to \{\beta : \beta\prec \alpha\}$ witnessing the $\alpha$-treeness of $T$ is obtained by letting $f(\sigma)$ be the least $\beta \prec \alpha$ such that  $\sigma\in g(\beta)$. 
\end{proof}

\begin{remark}\rm
In $\atro + \mathit{WO}(\alpha)$ every $\beta$-tree for $\beta \preceq \alpha$ is suitable.
\end{remark}

The next task, which presents itself, is to translate the formulas of $\ls$ that belong to the collection $\Bcal$ into the language of second order arithmetic, $\ltwo$.

\begin{definition}\rm \label{d:Bcalltwo}
For convenience assume that we have injections $x\mapsto T_x$ and $X\mapsto T_X$  
from set variables $x$ and relation variables $X$ of $\mathcal{L}^*$ to second order variables of $\ltwo$ in such a way that $T_x$ is always different from any $T_X$. The latter provides a translation of the terms of $\mathcal{L}^*$ except for the constants $\underline{\emptyset}$ and $\underline{\omega}$. These can just be translated as  $\{\langle\rangle\}$ and $\omega^*$, respectively (see Definition~\ref{MR6}). 

Since the relation $\in^*$ of Definition~\ref{MR3} has only been defined for trees, let us agree that henceforth a formula $S\in^* T$ will be considered a shorthand for 
$S,T\mbox{ trees } \lland S\in^* T$. 

The translation $^*$ from $\Bcal$ to $\ltwo$ is then effected as follows.
\begin{enumerate}
\item 
$(a\in \underline{\omega})^*$ is $T_a\in^*\omega^*$. $(\neg a\in \underline{\omega})^*$ is the negation of $(a\in \underline{\omega})^*$. $(a\in \underline{\emptyset})^*$ is $0=1$ while $(\neg a\in \underline{\emptyset})^*$ is $0=0$. 
\item
$(a\in x)^*$ is $T_a\in^* T_x$.   $(a \notin x)^*$ is the negation of $(a \in x)^*$.
\item
$(U(a))^*$ is $T_a\in^* T_U\,\wedge\, T_U\subseteq^*\omega^*$.  $(\neg U(a))^*$ is the negation of $(U(a))^*$.
\item
$(M_{\alpha}(a))^*$ is {\em $T_a$ is an $\alpha$-tree}.  $(\neg M_{\alpha}(a))^*$ is {\em $T_a$ is not an $\alpha$-tree}.
\item
$(F\lor G)^*$ and  $(F \land G)^*$ are  $F^* \lor G^*$ and  $F^* \land G^*$, respectively.
\item
$((\exists x \in a)F[x])^*$ and $((\forall  x\in a)F[x])^*$ are
\[ \exists T_x(T_x\in^* T_a \lland F^*[T_x]) \quad\mbox{and}\quad 
\forall T_x(T_x\in^* T_a \tto F^*[T_x]) \]
respectively, where $ F^*[T_x]$ stands for $(F[x])^*$; note that $^*$ replaces $x$ by $T_x$. 
\item 
$(\exists x^{\alpha}F[x])^*$  and $(\forall x^{\alpha}F[x])^*$ are
\[ \exists T_x\,(T_x\mbox{ $\alpha$-tree} \lland F^*[T_x]) \quad\mbox{and}\quad
  \forall  T_x\,(T_x\mbox{ $\alpha$-tree} \tto F^*[T_x]), \]
respectively.
\item 
$(\exists X F[X])^*$ and $(\forall X F[X])^*$ are
\begin{gather*}
\exists T_X(T_X\mbox{ $\omega{+}1$-tree} \lland T_X\subseteq^*\omega^* \lland F^*[T_X])
\end{gather*}
and
\begin{gather*}
\forall  T_X(T_X\mbox{ $\omega{+}1$-tree} \lland T_X\subseteq^*\omega^*\tto  F^*[T_X]),
\end{gather*} 
respectively. Here $ F^*[T_X]$ stands for $(F[X])^*$, noting that $^*$ replaces $X$ by $T_X$.
\end{enumerate}
\end{definition}

In order to show that the collapsed  derivations of Theorem~\ref{t:refined-collapsing} prove true statements when subjected to the interpretation of the previous, we need to show that the rule $\BC$ preserves truth. This will mainly be a consequence of
$\Sigma^1_1$-$\mathsf{AC}$ .

\begin{lemma}\label{MR9} 
Let $F[x_0,\ldots ,x_q,U_0,\ldots,U_r]$ be a $\Sigma$-formula of $\mathcal{L}_2^{set}$ with all free variables indicated,  and $\ell_F$ be the length of the latter formula as a string of symbols. 
Let
\[ F^{\alpha}[x_0,\ldots ,x_q,U_0,\ldots,U_r] \quad\mbox{and}\quad 
   F^{(a)}[x_0,\ldots ,x_q,U_0,\ldots,U_r] \]
be the results of replacing the unbounded existential set quantifiers $\exists x$ in the formula by $\exists x^{\alpha}$ and $\exists x\in a$, respectively.  

Arguing in the theory $\TTa$,  we have that whenever  $T_0,\ldots,T_q$ are $\alpha$-trees and $Q_0,\ldots ,Q_r$ are trees such that $Q_i\subseteq^* \omega^*$ (so actually $\omega{+}1$-trees), and the $*$-translation of
\[ F^{\alpha}[x_0,\ldots ,x_q,U_0,\ldots,U_r] \]
with $x_i$  and $U_j$ replaced by $T_i$ and $Q_j$, respectively,  is true, then the $*$-translation of 
\[ \exists y^{\alpha+\ell_F}F^{(y)}[x_0,\ldots ,x_q,U_0,\ldots,U_r],\]
again  with $x_i$ replaced by $T_i$ and $U_j$ replaced by $Q_j$, holds true as well.
\end{lemma}

\begin{proof}\setcounter{equation}{0}
We proceed by (meta)  induction on $\ell_F$. The most interesting case arises when this formula starts with a bounded universal quantifier, i.e., if it is of the form
\[ \forall v\in x_0\,F^{\alpha}_0[x_0,\ldots ,x_q,v,U_0,\ldots,U_r]. \]
 So assume that $T_0,\ldots,T_r$ are $\alpha$-trees and $Q_0,\ldots ,Q_r$ are trees such that $Q_i\subseteq^* \omega^*$ and
 \[ (\forall \langle n\rangle\in T_0)
    G^{\alpha}_0[T_0,\ldots ,T_q,T_0^{\langle n\rangle},Q_0,\ldots,Q_r] \]
 holds, where $G^{\alpha}_0$ is the $*$-translation of  $F^{\alpha}_0$. Inductively we then have %
\begin{gather}\label{MRc}
(\forall \langle n\rangle\in T_0) \exists S(S \mbox{ an $\alpha{+}\ell_{F_0}$-tree} \;\land \;
G^{(S)}_0[T_0,\ldots ,T_q,T_0^{\langle n\rangle},Q_0,\ldots,Q_r]),
\end{gather}
where $G^{(S)}_0[T_0,\ldots ,T_q,T_0^{\langle n\rangle},Q_0,\ldots,Q_r]$ results from $G^{\alpha}_0[T_0,\ldots ,T_q,T_0^{\langle n\rangle},Q_0,\ldots,Q_r]$ by replacing any part
of the form $\exists Q[Q\mbox{  an $\alpha$-tree} \lland \ldots Q\ldots]$ by $(\exists \langle n\rangle \in \mathcal{S})(\ldots \mathcal{S}^{\langle n\rangle}\ldots)$. 

Now it follows from Lemma~\ref{MRalpha} that the part $(\ldots)$ in (\ref{MRc}) is $\Delta^1_1$. As a result, we may apply $\Sigma^1_1$-$\mathsf{AC}$ (which is  provable in $\Pi^1_1\mbox{-}\mathsf{CA}_0$) to conclude that there exists a set $Z$ such that for all $\langle n\rangle \in T_0$,
\begin{gather}\label{MRa} 
Z_{(n)}\mbox{ an $\alpha{+}\ell_{F_0}$-tree} \;\land\; 
G^{(Z_{(n)})}_0[T_0,\ldots ,T_q,T_i^{\langle n\rangle},Q_0,\ldots,Q_r]
\end{gather}
where $Z_{(n)}\,=\,\{k : (n,k)\in Z\}$. 
 
Using $\Sigma^1_1$-$\mathsf{AC}$  again, we can also single out a sequence of functions
\[ f_{2^n3^k} : Z_{(n)}^{\langle k\rangle} \tto \alpha{+}\ell_{F_0} \]
 such that $f_{2^n3^k}$ witnesses that $Z_{(n)}^{\langle k\rangle} $ is an
 $\alpha{+}\ell_{F_0}$-tree whenever  $\langle n\rangle \in T_0$ and 
 $\langle k\rangle \in Z_{(n)}$. 
   
Now define a tree $\mathcal S$ by 
\[ \mathcal{S} \;:=\; \{ \langle\rangle \} \;\cup\;
   \{\langle 2^n3^k \rangle \star \sigma : \langle n\rangle \in T_0 \lland  \langle k\rangle\in Z_{(n)}
   \lland \sigma \in Z_{(n)}^{\langle k\rangle}\}. \]
By design, $S$ is a tree and, moreover,  $S$ is an $\alpha{+}\ell_F$-tree as witnessed by the function $f\colon S\to \alpha{+}\ell_F$ defined by
\[ f(\langle\rangle) \;:=\; \alpha{+}\ell_{F_0} \quad\mbox{and}\quad
   f(\langle 2^n3^k \rangle \star \sigma) \;:=\; f_{2^n3^k}(\sigma). \]
 As a consequence of (\ref{MRa}) and the fact that $F[x_0,\ldots ,x_q,U_0,\ldots,U_r]$ is a $\Sigma$-formula, we infer that
\[ (\forall \langle n\rangle \in T_0)
   G^{\mathcal S}_0[T_0,\ldots ,T_q,T_0^{\langle n\rangle},Q_0,\ldots,Q_r] \]
and thus, as desired,  the truth of the $*$-translation  of $\exists y^{\alpha{+}\ell_F}F^{(y)}(x_0,\ldots ,x_q,U_0,\ldots,U_r)$ with the substitutions $x_i\mapsto T_i$ and $U_j\mapsto Q_j$ follows.
\end{proof}

\subsection{Conservativity}

It remains to ascertain that $\kp + \poocas$ is conservative over $\TT$ for formulas in the language of second order arithmetic. More precisely, if $\Fcal$ is a formula of second order arithmetic (as for instance defined in \cite[I.2]{simpson09}) and $\Fcal_0$ denotes its natural translation into the language $\mathcal{L}_2^{set}$, then we aim to show that
\[ \kp + \poocas \vvdash \Fcal_0 \impl \TT \vvdash \Fcal. \]

Now assume  $\kp + \poocas \vdash \Fcal_0$.  The latter being a finite deduction, it follows from what we have established in the previous sections that we can determine  fixed naturals $n,p$ such that 
$\Hcal \prov{\beta}{\eta,p} \Fcal_0$ for some  derivation operator $\Hcal$ and $\beta,\eta\in C(\omega_n[\Omega+1],0)$. This follows from Theorems~\ref{t:refined-embedding},
\ref{t:refined-pred-cut-el}, and \ref{t:refined-collapsing}. To dilate on this, note that because of the finiteness of the deduction in $\kp + \poocas$ we can find this a priori bound 
$\omega_n[\Omega+1]$ so that  the entire ordinal analysis, commencing with the embedding of this finite deduction into the infinitary proof system, solely uses ordinals from $C(\omega_n[\Omega+1],0)$.

Moreover, we can carry this ordinal analysis out in the background theory $\TT$ as there exists an order preserving map from $C(\omega_n[\Omega+1],0)$ into the segment of ordinals  $<\psi(\omega_{n+1}[\Omega+1])$ for which transfinite induction is available in this theory. Such a map exists,  for if $\delta,\eta\in C(\omega_n[\Omega+1],0)$ with $\delta<\eta$ one has
\[ \psi(\omega_n[\Omega+1]+\omega^{\Omega+\delta}) < \psi(\omega_n[\Omega+1]+\omega^{\Omega+\eta}) < \psi(\omega_{n+1}[\Omega+1]) \]
by Lemma~\ref{l:o3}(2).

In light of the previous, the next step will be to show that $\Hcal \prov{\beta}{\eta,p} \Fcal_0$ entails that $\Fcal$ is true, all the while working in our background theory $\TT$. This is where the up to now neglected parameter $p$ comes into its own in that the idea is to employ a formal truth predicate for formulas of length $<p$. However, it is not possible to just focus on formulas $\Fcal_0$, where $\Fcal$ resides in the language of second order arithmetic, since in showing that $\Hcal \prov{\beta}{\eta,p} \Fcal_0$ implies the truth of $\Fcal$ we shall induct on $\beta$ and the pertaining derivation is usually not cut-free, so we have to take formulas from $\Bcal$ of length $<p$ into account.

It is perhaps noteworthy that if  $\Hcal \prov{\beta}{\eta,p} \Gamma$ with $\Gamma$ a set of formulas in $\Bcal$ and $\beta,\eta<\Omega$, then all formulas occurring in the derivation 
must be in $\Bcal$, too, as there can be no cuts in it involving formulas with unbounded set quantifiers since their ranks would be $\geq \Omega$ (see Definition~\ref{d:rs-formulas}  items (7) and (8)). 

The translation of the formulas from $\Bcal$ into $\ltwo$ formulas has been introduced in Definition~\ref{d:Bcalltwo}. The purpose of the following definition is to fix an arithmetized truth definition for $\Bcal$-formulas of length $< p$ in $\TT$. 

\begin{definition}\rm \label{W-Definition}
Let $Z$ be a set of naturals such $(Z)_k$ is  an $\alpha$-tree for all $k$. 
We'd like to engineer a formula $\truth_p(x,X)$ of $\ltwo$ such that for all $\Bcal$-formulas $F[x_1,\ldots,x_r,U_1,\dots,U_s] $  of length $<p$ with all free variables exhibited, 
\begin{gather}
\truth_p(\goed{F[x_1,\ldots,x_r,U_1,\dots,U_s]},Z) \;\lra\; \tilde{F}[Z_{(1)},\ldots Z_{(r)},Z_{(r+1)},\ldots, Z_{(r+s)}] \tag{*}
\end{gather} 
where $\goed{F[x_1,\ldots,x_r,U_1,\dots,U_s]}$ stands for the G\"odel number of $F[x_1,\ldots,x_r,U_1,\dots,U_s]$ whilst
\[ \tilde{F}[Z_{(1)},\ldots Z_{(r)},Z_{(r+1)}, \ldots, Z_{(r+s)}] \]
denotes the formula obtained from  $(F[x_1,\ldots,x_r,U_1,\dots,U_s])^*$ by replacing $T_{x_i}$ by $(Z)_i$ and $T_{U_j}$ by $(Z)_{r+j}$. Thus the sections $(Z)_k$ of $Z$ furnish an assignment of $\alpha$-trees to the free variables of $(F[x_1,\ldots,x_r,U_1,\dots,U_s])^*$. 
\end{definition}

The formal definition of such a truth predicate is a standard but cumbersome procedure. A place in the literature, where one finds this carried out in detail, is Takeuti \cite[CH.~3,19]{takeuti87}, and another is Troelstra \cite[1.5.4]{troelstra73}.

\begin{theorem}\label{Wahrheit}  
Fix an $n$ and $p$.  Let $\mathit{TV}(Z,\alpha)$ be short for ``$\forall k\,(Z)_k\mbox{ is an $\alpha$-tree}$''. Then  $\TT$ proves that for all sequents of formulas $\Gamma$ in $\Bcal$ consisting of formulas of length $<p$ and $\alpha,\beta,\rho\in C(\omega_n[\Omega+1],0)\cap\Omega$ that
$\Hcal \prov{\beta}{\eta,p} \Gamma$ implies
\[ \forall Z(\mathsf{TV}(Z,\alpha)\to  \truth_p(\goed{\bigvee \Gamma},Z)) \]
where $\bigvee \Gamma$ stands for the disjunction of  all formulas of $\Gamma$, but if
$\Gamma$ is empty let $\bigvee\Gamma$ be $\underline{\emptyset}\in \underline{\emptyset}$. \end{theorem}

\begin{proof}
 We reason in $\TT$ by induction on $\beta$. The axiom cases are obvious. If 
 $\Hcal \prov{\beta}{\eta,p} \Gamma$ is the result of an inference of a form other than $\BC$, then this follows immediately from the induction hypothesis applied to the immediate subderivations. Note also that in case of a cut, the cut formulas belong to $\Bcal$ and have lengths $<p$. Observe also that the derivation cannot contain $\sref$ inferences since
 $\beta<\Omega$. 
 
 So it remains to deal with the case where the last inference is an instance of $\BC$. Fortunately, this is what Lemma~\ref{MR9} is really about. The latter shows that if the premise $\Theta$ of an instance of $\BC$ is true under an assignment $Z$, then so is the conclusion.
 \end{proof}
 
 \begin{theorem}\label{Trompete} $\kp + \poocas$ is conservative over  $\TT$ for formulas of second order arithmetic. More precisely, if  $\Fcal$ is sentence of second order arithmetic (i.e. $\ltwo$) and  $\Fcal_0$ denotes its natural translation into the language $\ltset$, then
\[ \kp + \poocas \vdash \Fcal_0 \impl \TT\vdash \Fcal. \]
\end{theorem}

\begin{proof}
Assume $\kp + \poocas \vdash \Fcal_0$. As elaborated on before, with the help of the results in subsection~\ref{ss:finite-bounds} and Theorem~\ref{Wahrheit},  it follows that
\[ \TT \vdash \forall Z(\mathsf{TV}(Z)\to \truth_p(\goed{ \Fcal_0},Z)). \]
Thus, in light of (*)  and noting that $\Fcal_0$ has no free variables, 
\[ \TT\vdash \Fcal_0^* \]
where $\Fcal_0^*$ is the translation of $\Fcal_0$ according to Definition~\ref{W-Definition}. It remains to establish the relationship between $\Fcal_0^*$ and $\Fcal$.

$\Fcal_0^*$ arises from $\Fcal$ by translating numerical quantifiers $Q n$ as ranging over the immediate subtrees of $\omega^*$ and second order quantifiers as ranging over $\omega{+}1$-trees $T$ such that $T\subseteq^*\omega^*$. Now, as the naturals with their  ordering are isomorphic to the immediate subtrees of $\omega^*$ ordered via $\in^*$ and also the collection of  sets of natural numbers is naturally isomorphic to the collection of $\omega{+}1$-trees $T$ such that $T\subseteq^*\omega^*$, it follows that $\Fcal_0^*$ implies $\Fcal$, completing the proof. Anyone insisting on a more formal proof is invited to proceed by induction on the buildup of $\Fcal$.
\end{proof}

\begin{remark}\rm
In this article we do not compute the proof-theoretic analysis of the system $\TT$. However, it seems that the ordinal analysis of $\TT$ simply has to follow that of  $\pooca$ with
$\varepsilon_0$ replaced by $\BH$. So we conjecture that it is the ordinal $\psi_0(\Omega_\omega \cdot \BH)$ in the terminology of Buchholz and Sch\"utte \cite{buchholz-schuette88}.

For the ordinal analysis of $\pooca$ see, for example, the relevant chapters in Buchholz, Feferman, Pohlers and Sieg \cite{bfps81}, Buchholz and Sch\"utte \cite{buchholz-schuette88}, and J\"ager \cite{j79}.
\end{remark}

\section{Extensions}

Thus far we have investigated  what happens if one adds Kripke-Platek set theory to the subtheory of second order arithmetic based on $\Pi^1_1$-comprehension. There obtains a certain analogy to what Barwise  \cite{barwise75} called the \emph{Admissible Cover}, $\mathbb{C}\mathrm{ov}_{\mathfrak M}$, of a basic structure $\mathfrak M$. In his case, the basic structures $\mathfrak M$ were models of set theory.  $\mathbb{C}\mathrm{ov}_{\mathfrak M}$ is the intersection of all admissible sets which cover $\mathfrak M$ (see \cite{barwise75}, Appendix 2.1).\footnote{The proper proof-theoretic counterpart of the admissible cover was developed in J\"ager \cite{j86a}.}

In our context, the admissible cover amounts to grafting the theory $\kp$ onto a subsystem of second order arithmetic, $T$. This could be called the {\em proof-theoretic admissible cover} of $T$. There is, however, a crucial difference between Barwise's model-theoretic construction and the proof-theoretic one employed in this paper. In the former the basic structure one starts from remains unchanged in a strong sense when building its admissible cover in that no new subsets of $\mathfrak M$ become available in $\mathbb{C}\mathrm{ov}_{\mathfrak M}$ (see \cite{barwise75} Appendix Corollary 2.4), whereas in the proof-theoretic case the axioms of the basic theory $T$ will interact with the axioms of the set theory $\kp$, witnessed by the fact that, in general,  more theorems of second order arithmetic become deducible in
$T+ \kp$ than in $T$ alone.

It is also interesting to investigate how the proof-theoretic admissible cover plays out in the case of other well-known subsystems of second order arithmetic. 
It turns out that a certain pattern emerges. Moreover, the techniques developed in this paper,  when combined with insights 
from the literature,  suffice to get these additional results.

\begin{definition}\rm
We will focus on two well-studied theories. 
\begin{enumerate}
\item
$\atro$ with its signature axiom of \emph{arithmetical transfinite recursion} is the fourth system of the ``big'' five of reverse mathematics (see \cite[V]{simpson09}). 
\item
The second system we will consider is  traditionally called the theory of  \emph{bar induction}, 
$\bi$, by proof theorists. In \cite[VII.2.14]{simpson09} it is denoted by $\Pi^1_{\infty}\mbox{-}\mathsf{TI}_0$. The axioms of $\mathsf{BI}$ are those of $\mathsf{ACA}_0$ plus the scheme of transfinite induction
\[ \forall X(\mathit{WO}(X) \to \mathit{TI}(X,A)) \]
for every $\ltwo$ formula $A[x]$, where $\mathit{TI}(X,A)$ stands for the formula
\[ \forall u((\forall v <_X u)A[v] \to A[u]) \tto \forall u A[u] \]
with $\mathit{WO}(X)$ expressing that the ordering $<_X$ defined by  $v<_X u:\Leftrightarrow 2^v\cdot3^u\in X$ is a well ordering.
\end{enumerate}
\end{definition}

$\atros$ and $\bis$ are the axiom schemas -- formulated in the language $\ltset$ -- that comprise the natural translations of all instances of arithmetical transfinite recursion and transfinite induction, respectively. If $\Ccal$ is a collection of sentences of $\ltwo$ and $T_1$, $T_2$ are theories of the language $\ltwo$ or $\ltset$, we write $T_1\equiv_\Ccal T_2$ to convey that $T_1$ and $T_2$ possess the same $\Ccal$-theorems (perhaps modulo the translation of $\Ccal$ into $\ltset$).

\begin{theorem}\label{THm1}
$\kp + \atros \;\equiv_{\Pi^1_{\infty}}\; \atro + \mathit{TI}({<}\BH)$.
\end{theorem}

\begin{proof}
The proof of Theorem \ref{Trompete} essentially carries over with $\atros$ replacing $\poocas$ since the crucial Lemma~\ref{MR9} is also provable in $\atro$ as the latter theory proves $\Sigma^1_1$-$\mathsf{AC}$, too  (this is an old result, see \cite[Theorem V.8.3]{simpson09}).
\end{proof}

\newpage

\begin{theorem}\label{THm2} \quad
\begin{enumerate}[(i)]
\item
$\kp + \bis \;\equiv_{\Pi^1_{\infty}}\; \bi$.
\item
$\kp + \bis \;\equiv_{\Pi^1_{1}}\; \kp$. 
\item
$\kp + \atros \;\equiv_{\Pi^1_{1}}\; \kp$. 
\end{enumerate}
\end{theorem}

\begin{proof}
(i) The proof of Theorem \ref{Trompete} also works  with $\bi$ in lieu of  $\pooca$ since
$\atro$ is a consequence of $\bi$ (see \cite[Corollary VII.2.19]{simpson09}. So we infer that $\kp + \bis \equiv_{\Pi^1_{\infty}} \bi+\mathit{TI}({<}\BH)$. However, it is well known from the proof-theoretic literature that $\bi$ already proves $\mathit{TI}({<}\BH)$, yielding (i) (for more details see \cite{f70a,bfps81,buchholz-schuette88}). 

(ii) By J\"ager \cite{j82a} and the papers cited in the previous line,  it is known that $\kp$ and $\bi$ share the same proof-theoretic ordinal. Moreover, from the ordinal analyses of these two theories it can be inferred that they prove the same $\Pi^1_1$-theorems. 

The degree of conservativity, though,  cannot be improved much beyond this level as
$\atro$ has a $\Pi^1_2$-axiomatization and $\kp$ does not prove $\mathsf{ATR}_0$. To see this, first note that $\mathsf{KP}$ has a model $\mathbb{H}\mathrm{YP}_{\mathfrak N}$ which is the intersection of all admissible sets that contain the standard structure $\mathfrak N$ of the naturals as a set (see \cite[Theorem 5.9]{barwise75} ). The subsets of $\Nbb$ in $\mathbb{H}\mathrm{YP}_{\mathfrak N}$ are the hyperarithmetical sets. Then we can proceed, for example, in one of the following two ways:
\begin{enumerate}
\item[-]
The collection of the hyperarithmetical sets, as apparently proved by Kreisel, constitutes the smallest $\omega$-model of $\Sigma^1_1$-$\mathsf{AC}$, but it cannot be a model of  $\atro$ as it would have to have a proper $\omega$-submodel that is again an $\omega$-model of $\atro$, and  thus of  $\Sigma^1_1$-$\mathsf{AC}$. The result about $\omega$-models of  $\mathsf{ATR}_0$ is due to Quinsey \cite[pages 93--96]{quinsey80}  (for more details see \cite{simpson09} Theorem VIII.6.12 and the Notes for $\S$VIII.6).
\item[-]
Alternatively, observe that the theory $\atro$ is equivalent to the theory $\fpo$ which extends
$\mathsf{ACA}_0$ by axioms stating that any positive arithmetic operator has a fixed point; see Avigad \cite{avigad96a}). But the hyperarithmetical sets do not constitute a model of $\fpo$ according to Probst \cite[II.2.4]{probst05a} and Gregoriades \cite{gregoriades19a}.
\end{enumerate}

(iii) By Theorem \ref{THm1} we have  $\atro^* + \kp  \equiv_{\Pi^1_{\infty}} \atro +\mathit{TI}({<}\BH)$. Since $\atro + \mathit{TI}({<}\BH)$ is a subtheory of $\bi$  and 
$\bi^* + \kp \equiv_{\Pi^1_{1}} \kp$, it follows that $\atro^* +\kp \equiv_{\Pi^1_{1}} \kp$.
\end{proof}

In connection with the previous results we would also like to mention Sato \cite{sato18a}. Among other things, he states there a special case of the previous Theorem~\ref{THm2}(ii), namely that the addition of a $\Pi^1_2$ theorem of $\bi$ to $\kp$ does not increase the consistency strength of the augmented theory.%
\footnote{What we call $\kp$ is called $\kp\omega$ in \cite{sato18a}.}
Further interesting work about the relationship between Kripke-Platek set theory and $\Pi^1_1$ comprehension on the natural numbers is presented in Simpson \cite{simpson18a}. There the interplay between $\kp$  and $\pooca$ is studied from a recursion- and model-theoretic perspective. However, it does not provide the exact proof-theoretic strength of $\kp + \poocas$.

\bigskip

\bibliography{kppiooca}
\bibliographystyle{amsplain}

\noindent\textbf{Addresses}
\medskip

\noindent
Gerhard J\"ager\\
Institut f\"ur Informatik\\
Universit\"at Bern\\
Switzerland \\
gerhard.jaeger@inf.unibe.ch

\medskip

\noindent
Michael Rathjen\\
Department of Pure Mathematics\\
University of Leeds\\
England \\
M.Rathjen@leeds.ac.uk
\end{document}